\documentclass[11pt,fullpage]{article}
\usepackage{multicol,wrapfig,amsmath,subfigure}
\usepackage{epstopdf}
\usepackage{multirow}
\usepackage{amsmath}
\usepackage{amsfonts,amssymb}
\usepackage{amsthm}
\usepackage{graphics,graphicx}

\usepackage[hyphens,spaces,obeyspaces]{url}

\usepackage{hyperref}
\usepackage[right = 2.5cm, left=2.5cm, top = 2.5cm, bottom =2.5cm]{geometry}
\usepackage{tikz-cd}
\usepackage{slashed}
\usepackage{enumerate}
\usepackage{comment}
\usepackage{authblk}
\usepackage{mathtools}
\newtheorem{prop}{Proposition}[section]
\newtheorem{lemma}{Lemma}[section]
\newtheorem{defn}{Definition}[section]
\newtheorem{thm}{Theorem}[section]
\newtheorem{cor}{Corollary}[section]

\newtheorem{question}{Question}
\newtheorem{exercise}{Exercise}
\newtheorem{example}{Example}[section]

\newcommand{\R}{\mathbb{R}}
\newcommand{\C}{\mathbb{C}}

\newcommand{\Delbar}{\overline{\partial}}
\newcommand{\Endpoint}{{\mathfrak{Ep}}}

\graphicspath{ {images/} }

\begin{document}

\title{Directed immersions for complex structures}
\author{Tobias Shin}

\begin{center}
\Large\textbf{Directed immersions for complex structures}
\linebreak
\text{\normalsize{Tobias Shin}}\\
\end{center}
\begin{abstract}
We analyze the differential relation corresponding to integrability of almost complex structures, reformulated as a directed immersion relation by Demailly and Gaussier. Combining results of Clemente \cite{clemente}, we show that applying $h$-principle techniques yields the following statement: for an almost complex manifold with arbitrary metric $(X, J, g)$, and for $\epsilon > 0$, there exists a smooth function $f: X \rightarrow \R$ and almost complex structure $J'$ on $X$ such that $J$ and $J'$ are $C^0$-close on the graph of $f$ with respect to the extended metric on $X \times \R$, and such that the Nijenhuis tensor of $J'$ on the graph has pointwise sup norm less than $C\epsilon$, where $C$ is a constant depending only on $J$ and $g$. This is an updated version of a previous preprint titled "Almost complex manifolds are (almost) complex". 
\end{abstract}
\section{Introduction}
In their paper \cite{demailly}, Demailly and Gaussier construct, for a given complex dimension $n$, a universal space $Z$ with an algebraic distribution $D$ for which all almost complex $n$-manifolds immerse into, transverse to $D$. More precisely, they prove the following theorem.
\begin{thm} (Demailly, Gaussier)
For all integers $n \geq 1$, there exists a complex affine algebraic manifold $Z$ of dimension $N = 38n^2 + 8n$ possessing an anti-holomorphic algebraic involution and an algebraic distribution $D \subset TZ$ of codimension $n$, for which every compact $n$-dimensional almost complex manifold $(X,J)$ admits an embedding $f: X \hookrightarrow Z^\R$ transverse to $D$ and contained in the real part of $Z$, such that $J = J_f$, where $J_f$ denotes the almost complex structure on $TZ/D$ pulled back under $f$.
\end{thm}
The space $Z$ they construct is a combination of Grassmannians and twistor spaces, built in such a way that essentially ``globalizes" the local picture relating Frobenius integrability with the Nijenhuis tensor, via Whitney embedding. Moreover, they give a criterion for when a given almost complex structure $J$ is integrable, with respect to this setup.
\begin{thm}\label{DemaillyTheorem} (Demailly, Gaussier)
For every compact $n$-dimensional \textnormal{integrable} complex manifold $(X,J)$, there exists an embedding $X \hookrightarrow Z^\R$ transverse to $D$, contained in the real part of $Z$, such that
\end{thm}
\begin{enumerate}[(i)]
\item \textit{$J = J_f$ and $\overline{\partial}_{J}f$ is injective};
\item Im$(\overline{\partial}_Jf)$ \textit{is contained in the isotropic locus $I$ of the torsion operator $\theta$ of $D$, the intrinsically defined algebraic locus in the Grassmannian bundle $Gr(n,D) \rightarrow Z$ of complex $n$-dimensional subspaces in $D$ consisting of those subspaces $S$ such that $\theta|_{S\times S} = 0.$}
\end{enumerate}
The torsion operator $\theta$ is defined as
\begin{equation*}
\theta(X,Y) = \lbrack X, Y\rbrack \ \textnormal{mod} \ D.
\end{equation*}
Note that by taking the quotient, we obtain a skew symmetric bilinear \textit{tensor}; we may view it as a holomorphic section of the bundle $\wedge^2D^{*} \otimes (TZ/D)$. The inclusion condition (ii) Im$(\overline{\partial}_Jf) \subset I$ is actually necessary and sufficient for integrability of $J_f$. Here, $\overline{\partial}_Jf = \frac{1}{2}(df + J_Z\circ df \circ J_f)$ where $J_Z$ is the fixed complex structure on $Z$. It was first suggested by Demailly at James Simons' 80th birthday conference to use the above theorems as a strategy to find obstructions to \textit{moving into} the isotropy locus; this would give topological obstructions to having an integrable complex structure.
\\ \\
In the language of Gromov and the $h$-principle \cite{gromov}, the above gives rise to a natural \textit{directed immersion} problem. As described in Eliashberg \& Mishachev \cite{eliashberg}, the setup is as follows: let $Gr_n(W)$ be the Grassmannian bundle of tangent $n$-planes to a manifold $W$ of dimension strictly larger than $n$ and $V$ be an $n$-dimensional manifold. Let $A \subset Gr_n(W)$ be an arbitrary subset. An immersion $f: V \rightarrow W$ is said to be an \textit{A-directed immersion} if the induced tangential lift $Gdf$ maps into $A$, where $Gdf$ sends $v$ to $df_v(T_vV)$. 
All of the above can similarly be done replacing the word \textit{immersion} with \textit{embedding}. 
\\ \\
In other words, the question of when an almost complex structure can be moved along a path of almost complex structures into one that is integrable is equivalent to a modified directed immersion problem with respect to $\Delbar$ of immersions into the universal space $Z$. One can then apply the philosophy of the $h$-principle and ask:
\begin{enumerate}[(i)]
\item Are there formal solutions to this differential relation?
\item Assuming (i), are there genuine solutions to this differential relation? Does it satisfy the $h$-principle?
\end{enumerate}
In fact, even if there were no obstructions to a formal solution, we already know that this relation \textit{fails} the $h$-principle: there exist almost complex manifolds in complex dimension 2 that have no integrable complex structures, as shown classically by Van de Ven \cite{van}. It is an open problem as to whether there exist such manifolds in higher dimensions.
\\ \\
This paper is organized as follows: we first give preliminary information about the differential relation and the corresponding subspace of the Grassmannian bundle. We prove in section \ref{formal} the following result.
\begin{thm}
There are always \textnormal{formal} solutions to the above differential relation. 
\end{thm}
In the previous version of this paper, the above theorem only held for complex dimensions up to 77. However, using results of Clemente \cite{clemente}, we can say it holds for all dimensions.
\\ \\
In section \ref{ample}, we prove that there is always a holonomic section into an arbitrarily small \textit{open neighborhood} of the relation, using the method of holonomic approximation and a microextension trick \cite{eliashberg}, assuming the existence of a formal solution. The previous version of this paper erroneously stated this meant there were almost complex structures with Nijenhuis tensors whose $C^0$ norms became arbitrarily small. First, a minor definition:

\begin{defn} Let $(X,g)$ be a Riemannian manifold and let $f: X \rightarrow \R$ be a smooth function. Let $\Gamma_f$ denote the graph of $f$ as a submanifold of $X \times \R$. Restrict the metric $g_{ij} \oplus dt^2$ on $X \times \R$ to the graph $\Gamma_f$ to obtain a metric $\hat{g}_f$ on $\Gamma_f$. We will call the metric $\hat{g}_f$ the \textbf{graph metric with respect to $f$}.  
\end{defn}  
The argument presented in the previous version of this paper then actually proves the following statement.
\begin{thm}\label{MainTheorem}
Let $(X, J, g)$ be an almost complex manifold with metric. For any $\epsilon > 0$, there exists a smooth function $f: X \rightarrow \R$ and almost complex structure $J'$ on $X$ such that
\begin{itemize}
\item $||\pi^*J - \pi^*J'||_{C^0} < K\epsilon$ where $\pi^*J$ and $\pi^*J'$ denote the pullbacks of $J$ and $J'$ on $\Gamma_f$ respectively, and where the norm is taken with respect to $\hat{g}_f$ and the constant $K$ only depends on $J$ and $g$
\item the Nijenhuis tensor of $J'$ on $\Gamma_f$ has $C^0$ norm less than $C\epsilon$, again with respect to $\hat{g}_f$, and where $C$ is a constant depending only on $J$ and $g$.
\end{itemize}
\end{thm}
In other words, given an almost complex structure $J$ and metric $g$, we can find a smooth function $f$ and a new almost complex structure $J'$ so that $J'$ approximates $J$ on the graph of $f$, and so that the Nijenhuis tensor of $J'$ has small norm on the graph of $f$.
\\ \\
\textbf{Acknowledgements.} The author would like to first and foremost thank his advisor Dennis Sullivan for constant encouragement and guidance throughout the work of this paper. The author would also like to thank Jean-Pierre Demailly, who gave a series of inspiring lectures at James Simons' 80th birthday conference at CUNY Graduate Center and whose ideas presented there strongly motivated this work; Andrew Sommese, who clarified his vanishing relative homotopy theorems in an exchange of helpful emails; and \'{A}lvaro del Pino, who explained the full details of proofs from \cite{eliashberg} concerning the $h$-principle in a very enlightening correspondence. This project is really a love letter to the Stony Brook math department family. The author would also like to thank Frederik Benirschke, Mark De Cataldo, Nathan Chen, Yoon-Joo Kim, Robert Lazarsfeld, Lisa Marquand, John Sheridan, Jason Starr, Sasha Viktorova, Ben Wu, and Ruijie Yang for teaching the author algebraic geometry; Michael Albanese, Jiahao Hu, and Aleksandar Milivojevic for teaching the author algebraic topology; Corey Bregman, Robert Bryant, Jae Ho Cho, Lisandra Hernandez-Vazquez, Ben McMillan, Jordan Rainone, and Dror Varolin for teaching the author complex and differential geometry; Dahye Cho, Mohamed El Alami, Yasha Eliashberg, Oleg Lazarev, Mark McLean, Tony Phillips, Ying Hong Tham, and Hang Yuan for teaching the author symplectic topology. Finally, the author thanks Scott Wilson, Luis Fernandez, and Alvaro del Pino Gomez for comments.

\section{The subspace $I$ and the relation $\mathcal{R}_I$}\label{prelim}
In this section we provide some preliminary information about the universal space $Z$, the distribution $D$, the subspace $I$, and the differential relation $\mathcal{R}_I$. For the entirety of this paper, we only discuss complex dimension.
\\ \\
We first recall the definition of jets and jet spaces, following Hirsch \cite{hirsch}.

\begin{defn}
Let $M, N$ be two smooth manifolds. An \textbf{r-jet} is an equivalence class $[x,f,U]_r$ of a triple $(x,f,U)$ where $U \subset M$ is open, $x \in U$, and $f:U \rightarrow N$ is a $C^r$ map. Two triples $(x,f,U)$ and $(x',f',U')$ are equivalent if $x = x'$ and if in some (and therefore any) pair of coordinate charts, $f$ and $f'$ have the same partial derivatives up to order $r$. We denote $[x,f,U]_r$ as $j^r_xf$, called the \textbf{r-th jet of f at x}. The space of all $r$-jets from $M$ to $N$ is denoted $J^r(M,N)$. 
\end{defn}
There is a natural projection map $J^{r+1}(M,N) \rightarrow J^r(M,N)$ which gives $J^{r+1}(M,N)$ the structure of an affine bundle over $J^r(M,N)$. The space $J^0(M,N)$ is simply the product $M \times N$.
\begin{defn} A \textbf{differential relation} $\mathcal{R}$ is a subset of $J^r(M,N)$.
\end{defn}
For this paper, we will mainly be concerned with the space $J^1(M,N)$. After choosing coordinates on $M$ and $N$, an element $j^1_x f$ of $J^1(M,N)$ can be represented as a tuple $(x,y,L)$ where $y$ is the image of $x$ under $f$ and $L \in $ Hom$_{\R}(T_xM, T_yN)$ is the derivative of $f$ at $x$. In other words, we can think of $L$ as a ``formal'' a priori derivative at $x$.
\\ \\
Let us illustrate some examples of differential relations below.
\begin{example} We write $\mathcal{R}_{imm} \subset J^1(M,N)$ to denote the differential relation of maps that are \textbf{immersions}. This relation consists of tuples $(x,y,L)$ such that $L$ is injective.
\end{example}

\begin{example} Let $N$ be equipped with a complex structure $J_N$. We write $\mathcal{R}_{real} \subset J^1(M,N)$ to denote the differential relation of maps that are \textbf{totally real}. This relation consists of tuples $(x,y,L)$ such that  $L_x(T_xX) \cap J_N(L_x(T_xX)) = 0$.
\end{example}

\begin{example} Let $N$ be equipped with a distribution $D$. We write $\mathcal{R}_{D} \subset J^1(M,N)$ to denote the differential relation of maps that are \textbf{transverse} to $D$. This relation consists of tuples $(x,y,L)$ such that the composition $T_xX \xrightarrow{L} T_yN \rightarrow T_yN/D_y$ is surjective.
\end{example}
Notice that the above three examples are subsets that are \textit{open} in the jet spaces. Before the next example, we introduce some notation.
\begin{defn} Let $Gr(r,s)$ denote the Grassmannian of $r$-planes in $s$-space. A point in $Gr(r,s)$ is an $r$-plane in a complex vector space of dimension $s$.
\end{defn}
For a vector space equipped with a complex structure $J$, we say an $r$-plane $S$ is complex if and only if it is $J$-invariant; that is, $J(S) \subset S$. All $r$-planes are complex $r$-planes in this paper, unless otherwise stated. 
\begin{defn} Let $M$ be an almost complex manifold. Let $Gr(n,TM)$ denote the \textbf{Grassmannian bundle} of $n$-planes in $TM$. For each point $x$ in $M$, the fiber of this bundle at $x$ is the Grassmannian $Gr(n,T_xM)$ of $n$-planes in the vector space $T_nM$. If $A$ is a subset of $Gr(n,TN)$, we write $A_x$ to denote the intersection of $A$ with the fiber over $x$.
\end{defn}
\begin{example} Let $A$ be a subset of $Gr(n,TN)$ and assume $M$ is of dimension $n$ and $N$ is an almost complex manifold of much larger dimension, with complex structure $J_N$. We write $\mathcal{R}_{A} \subset J^1(M,N)$ to denote the differential relation of maps that are $A$-\textbf{directed}. This relation consists of tuples $(x,y,L)$ such that $L_x(T_xM) \in A_y$. Notice this only makes sense if $L$ is injective and, as we require our planes to be complex, if $L_x(T_xM)$ is $J_N$-invariant.
\end{example}

\begin{defn} A smooth distribution $D$ of constant rank $k$ in $TM$ is a smooth section of the Grassmannian bundle $Gr(k, TM)$ over $M$.
\end{defn}
In this paper, all distributions will be of constant rank; distributions then also correspond to subbundles of $TM$. With this point of view, given a distribution $D$, it makes sense to discuss general $A$-directed mappings where $A$ is a subset of $Gr(\cdot, D)$.
\\ \\
Recall that for real analytic complex structures, there is a deep theorem first proven by Gauss relating integrability of complex structures to integrability of distributions in the sense of Frobenius.

\begin{thm}(Gauss-Frobenius) Let $U \subset \R^{2n}$ be an open set containing the origin and let $J$ be a real analytic complex structure on $U$. Let $U^\C$ denote the complexification $U + iU \subset \C^{2n}$ of $U$ and let $J^\C$ denote the complexification of $J$. Then $J^\C$ extends holomorphically to a small ball in $U^\C$. Moreover, $J$ is integrable if and only if the $+i$-eigenspaces of $J^\C$ on $U^\C$ are integrable in the sense of Frobenius.
\end{thm}
\textit{Sketch of proof}. We have that $J$ is a real analytic map $U \rightarrow $ Hom$(TU, TU)$. Thus it has a power series expansion in a neighborhood of the origin in $U$. Replacing real variables by complex variables, we obtain $J^\C$ which is a map $U \rightarrow $ Hom$(TU \otimes_\R \C, TU \otimes_\R \C)$, with the same power series as $J$. This power series converges in a small ball in $U^\C$ around the origin in $\C^{2n}$, since the series for $J$ is also convergent. Therefore we can extend $J^\C$ holomorphically to this ball in $U^\C$. This extension is also a complex structure and so has $+i$-eigenspaces defined on this ball around the origin.
\\ \\
To prove the last statement, consider $U$ embedded in $U + iU$ by the diagonal mapping. This mapping is totally real; moreover it is transverse to the distribution of $+i$-eigenspaces defined on the ball around the origin. If the $+i$-eigenspaces are integrable in the sense of Frobenius, we can find a holomorphic foliation of submanifolds that are tangent to the $+i$-eigenspaces and transverse to the diagonal. The leaves of this foliation then project to the holomorphic charts for which $J$ is constant. \qed
\\ \\
We now discuss the following deep theorem of Demailly and Gaussier \cite{demailly}.

\begin{thm} (Demailly, Gaussier)
For all integers $n \geq 1$, there exists a complex affine algebraic manifold $Z$ of dimension $N = 38n^2 + 8n$ possessing an anti-holomorphic algebraic involution and an algebraic distribution $D \subset TZ$ of codimension $n$, for which every compact $n$-dimensional almost complex manifold $(X,J)$ admits an embedding $f: X \hookrightarrow Z^\R$ transverse to $D$ and contained in the real part of $Z$, such that $J = J_f$, where $J_f$ denotes the almost complex structure on $TZ/D$ pulled back under $f$.
\end{thm}
\textit{Sketch of proof.} The space $Z$ is defined as the set of tuples $(z,S',S'',\Sigma',\Sigma'')$ in $\C^{8n}\times Gr(3n,8n)\times Gr(3n,8n)\times Gr(4n,8n)\times Gr(4n,8n)$ where $S'\subset\Sigma'$, $S''\subset\Sigma''$, and $\Sigma'\oplus\Sigma'' = \C^{8n}$. The space $Z$ is then a quasiprojective subvariety of this product of Grassmannians; it is a smooth variety of dimension $N = 38n^2 + 8n$. The flag decompositions are equivalent to a choice of complex structure on $\C^{8n}$ and $S'\oplus S''$ where $S' \subset\Sigma'$ and $S'' \subset\Sigma''$ correspond to the $+i$ and $-i$ eigenspaces respectively. The distribution $D$ is defined at a point $P = (z,S',S'',\Sigma',\Sigma'')$ as the set of tangent vectors $(v,u',u'',w',w'')$ where $v \in S'\oplus\Sigma''$ tautologically, with no other conditions on the other components. This gives $T_PZ/D_P \cong \Sigma'/S'$, so $D$ is corank $n$.
\\ \\
The embedding of an almost complex manifold $(X,J)$ into $Z$ is constructed by first embedding $X$ into the diagonal of $X\times X$, and then embedding $X\times X$ by the product embedding into euclidean space via Whitney. One then complexifies the normal bundle of this embedding and its complex structure $\widetilde{J}$, defined as the direct sum of the complex structures $J, -J,$ and $J_\Delta$ where $-J$ is the conjugate complex structure and $J_\Delta$ is the complex structure on the double normal bundle. The embedding $f$ into $Z$ is then by taking as the four subspaces the $+i$ and $-i$ eigenspaces of $\widetilde{J}$ on $\C^{8n}$ and the complexified normal bundle, in such a way that $\Sigma'/S' \cong T^{1,0}X$. Thus $Z$ and $D$ are defined tautologically so that $X$ embeds transverse to $D$, and so that the pullback complex structure agrees with the initial complex structure, by way of a real (which is a posteriori complex) isomorphism $f^*(TZ/D) \cong T^{1,0}X$. \qed
\\ \\
\textbf{Remark.} In \cite{demailly}, $Z$ is actually taken to be an open affine subset of the above space by removal of an appropriate subvariety, but it is shown that any almost complex manifold $X$ embeds into the above quasiprojective variety anyway.
\\ \\
The above construction is strongly motivated by the proof of Gauss-Frobenius; indeed, it is a globalization of the local picture in Gauss-Frobenius made universal by Whitney embedding. For full proof and construction, see \cite{demailly}. We will often call $Z$ the \textbf{universal space} and $D$ the \textbf{universal distribution} in this paper.
\\ \\
We now prove some properties about the universal space $Z$ and its universal distribution $D$. Namely, $Z$ is homogeneous.
\begin{lemma}\label{affineaction} There is a natural transitive group action of $\C^{8n}\times GL(8n,\C)$ on $Z$. There is another transitive group action, by the affine group  \textnormal{Aff}$(\C^{8n})$, on $Z$ such that $D$ is invariant; that is, for an affine element $g \in \textnormal{Aff}(\C^{8n})$, the differential of the multiplication map $\mu_g$ maps $D$ to $D$.
\end{lemma}
\begin{proof} Define the action of $\C^{8n} \times$ GL$(8n,\C)$ by $(b, A) \cdot (z, S', S'', \Sigma', \Sigma'') = (z + b, AS', AS'', A\Sigma', A\Sigma'')$. This is transitive since the action is transitive on all components. Define the action of Aff$(\C^{8n})$ by $(b, A) \cdot (z, S', S'', \Sigma', \Sigma'') = (Az + b, AS', AS'', A\Sigma', A\Sigma'')$. This is again transitive on all components. Moreover, for a fixed element $g = (b,A)$, we have a smooth map $\mu_g: Z \rightarrow Z$ which is multiplication by $g$. Recall that $D$ at a point $P = (z,S',S'',\Sigma',\Sigma'')$ is the set of $(v,u',u'',w',w'')$ where $v \in S'\oplus\Sigma''$ tautologically. Then, the map $\mu_g$ sends $(v,u',u'',w',w'')$ in $D_P$ to $(Av, \widetilde{u'}, \widetilde{u''}, \widetilde{w'}, \widetilde{w''})$ in $D_{\mu_g(P)}$, for some $\widetilde{u'}, \widetilde{u''}, \widetilde{w'}, \widetilde{w''}$. Thus $D$ is invariant under this action.
\end{proof}

Before we further analyze the structure of the algebraic distribution $D$ on $Z$, we must prove some lemmas regarding \textit{stably almost complex} manifolds.

\begin{defn} A smooth manifold $M$ is \textbf{stably almost complex} if $TM\oplus\epsilon_{\R}^k$ is a complex vector bundle for some $k$, where $\epsilon_{\R}^k$ is the trivial rank $k$ real vector bundle.
\end{defn}

We use the following lemma (2.1 in \cite{goertsches}).
\begin{lemma}\label{SACconnectsum} (Goertsches, Konstantis) Let $M$ and $N$ be two real $d$-dimensional manifolds, $M\#N$ their connect sum, and $p: M\#N \rightarrow M$ and $q: M\#N \rightarrow N$ the two collapsing maps. Then $T(M\#N)\oplus \epsilon_{\R}^d \simeq p^{*}TM\oplus q^{*}TN$. This implies the connect sum of two stably almost complex manifolds is stably almost complex.
\end{lemma}

The following observation is due to a discussion with M. Albanese and A. Milivojevic. We write $\#^{k} M$ to denote the $k$-fold connect sum of $M$. That is, $\#^{k} M = M \# ... \# M$ with $k$ copies of $M$.
\begin{lemma}\label{oddDimLineSubbundle} In every odd dimension, there exists examples of almost complex manifolds whose tangent bundles admit no complex line subbundles.
\end{lemma}
\begin{proof} We have that $S^5\times S^5$ is complex \cite{calabi} and so its connect sums are stably almost complex, by lemma \ref{SACconnectsum}. By Theorem 2 in \cite{yang}, we have that $\#^{25} S^5\times S^5$ is almost complex. Moreover, it is $2$-connected with nonzero Euler characteristic, and so it has no complex line subbundle in its tangent bundle. Since $S^6$ is also 2-connected with nonzero Euler characteristic, we have that products of $S^6$ with $\#^{25} S^5\times S^5$ are almost complex manifolds that do not admit any complex line subbundle in their tangent bundle. This gives every dimension of the form $3a + 5b$, which gives every odd dimension except 7. Finally, we have that $\mathbb{H}\mathbb{P}^2 \# \mathbb{H}\mathbb{P}^2 \# (S^4 \times S^4)$ is a 2-connected almost complex manifold with nonzero Euler characteristic in dimension 4 (Prop. 6 in \cite{muller}). Again by taking a product with $S^6$, we obtain dimension 7.
\end{proof}
\textbf{Remark.} The above argument only uses almost complex and stably almost complex manifolds. An earlier version of this paper asked a question: Does every odd dimensional \textit{complex} manifold admit a complex line subbundle in its tangent bundle? M. El Alami provided the author with a reference to a paper by Lu and Tian \cite{lutian}, wherein they prove that $\#^{25} S^3 \times S^3$ is a complex manifold. The answer to the question is then no. Lu and Tian in fact prove much more: they prove \textit{all} connect sums of $S^3 \times S^3$ with itself are complex. 

\begin{lemma}\label{nosubbundles} Consider the quotient bundle $TZ/D$ over $Z$. The quotient bundle $TZ/D$ has no proper complex subbundles.
\end{lemma}
\begin{proof} By universality of the construction of $D$ and $Z$ in \cite{demailly}, we have that for given dimension $n$, every almost complex $n$-manifold $X$ embeds into $Z$ by some map $f$ transverse to $D$ such that $f^*(TZ/D) \cong TX$. Therefore, to show that $TZ/D$ has no proper complex subbundles, it suffices to exhibit in each dimension $n$ an example of an almost complex manifold whose tangent bundle does not split into complex subbundles. However it is known that $\C\mathbb{P}^{2k}$ admits no complex subbundles and that $\C\mathbb{P}^{2k+1}$ only admits a complex line subbundle and its complement for all $k$ \cite{glover}. But by lemma \ref{oddDimLineSubbundle}, we have shown that every odd dimension has examples of almost complex manifolds whose tangent bundles have no complex line subbundles.
\end{proof}
Before we proceed, we recall the following basic lemma with its proof \cite{varolin}.
\begin{lemma}\label{commutingvf} Let $M$ be a complex manifold with complex structure $J$, and let $X$ and $Y$ be holomorphic vector fields. Then $J\lbrack X, Y \rbrack = \lbrack JX, Y\rbrack$.
\end{lemma}
\begin{proof} Recall that the Lie bracket of two vector fields is equivalent to taking the Lie derivative of one with respect to the other. Also using product rule, we have
\begin{equation*}
\begin{split}
\lbrack JX, Y\rbrack
&= -\mathcal{L}_{Y}(JX) \\ 
&= -J\mathcal{L}_{Y}X + (\mathcal{L}_{Y}J)X \\
&= J\lbrack X, Y\rbrack
\end{split}
\end{equation*}
with $\mathcal{L}_{Y}J$ vanishing since holomorphic vector fields admit flows $\phi_t$ that are holomorphic. Thus the flow associated to $Y$ therefore satisfies $\phi_t^*J = d\phi_{-t}\circ J\circ d\phi_t = J$. This gives that $\mathcal{L}_{Y}J = \frac{d}{dt}|_{t=0}  \phi_t^*J = \frac{d}{dt}|_{t=0}  J = 0$. 
\end{proof} 

We now come to the essential property of the distribution $D$ on $Z$.
\begin{defn} Let $D$ be a distribution on a manifold $M$. We say it is \textbf{bracket-generating} if the pointwise evaluations of iterated Lie brackets of vector fields tangent to $D$ eventually span $TM$.
\end{defn}

\begin{prop}\label{Disbracketgenerating}
The algebraic distribution $D$ on the universal space $Z$ is bracket-generating.
\end{prop}
\begin{proof} For a point $z$, consider the subspace $\widehat{D}_{z}$ defined as the span (over the ring of smooth $\R$-valued functions) of all iterated Lie brackets of vector fields tangent to $D$ evaluated at $z$. For any other $z' \in Z$, there exists $g \in $ Aff$(\C^{8n})$ such that $g\cdot z = z'$ by transitivity of the affine group action on $Z$, by lemma \ref{affineaction}. Since $D$ is invariant under this action, and biholomorphisms push forward Lie brackets, we have that $\widehat{D}_{z'} = dg_z(\widehat{D}_{z})$ independent of the choice of $g$. Thus $\widehat{D}$ is a well defined subbundle with $D \subset \widehat{D} \subset TZ$. To show $D$ is bracket-generating, we will prove that $\widehat{D} = TZ$. Notice immediately that $\widehat{D}$ cannot equal $D$, as then $D$ would be integrable in the sense of Frobenius. This would imply every almost complex \textit{structure} would be integrable.
\\ \\
Moreover $\widehat{D}$ is $J_Z$-invariant, since $D$ is holomorphic: first choose a local holomorphic frame for $D$. If $v \in \widehat{D}$, then we can write $v$ as the $\R$-linear span of evaluations of \textit{holomorphic} vector fields spanning $D$ and their iterated brackets. Then $J_Zv$ is in $\widehat{D}$ by lemma \ref{commutingvf}, since $J_Z$ commutes with all brackets and since $D$ is $J_Z$-invariant. We conclude $\widehat{D}$ is a complex subbundle in $TZ$.
\\ \\
If $\widehat{D}$ is a proper subbundle of $TZ$, then we obtain a proper quotient bundle $TZ/D \rightarrow TZ/\widehat{D}$. After choosing a hermitian metric, this gives us a proper complex subbundle of $TZ/D$. However this contradicts lemma \ref{nosubbundles}.
\end{proof}

\begin{defn} Let $D$ be a distribution on a smooth manifold $M$. The mapping $\theta: D \times D \rightarrow TM / D$ defined by
\begin{equation*}
\theta(X,Y) = \lbrack X, Y\rbrack \ \textnormal{mod} \ D.
\end{equation*}
where $X$ and $Y$ are vectors extended locally to vector fields, is the \textbf{torsion tensor} of $D$. This map is also called the \textbf{curvature} of $D$. 
\end{defn}
Notice that $\theta$ is exactly the obstruction for a distribution to be integrable in the sense of Frobenius and that it is a well defined tensor due to the quotient. In our specific case with the universal distribution $D$ on $Z$, the torsion tensor $\theta$ is algebraic.

\begin{defn} Let $Gr(n,D)$ be the Grassmannian bundle of $n$-planes over $Z$, with $D$ the universal distribution. Let $I$ denote the set of $n$-planes $S$ in $Gr(n,D)$ where $\theta|_{S \times S}$ is identically zero. We call $I$ the \textbf{isotropy locus} or \textbf{integrability locus}. We will also write $I_z = I \cap Gr(n,D_z)$ as the intersection of $I$ with the fiber of $Gr(n,D)$ over $z \in Z$. 
\end{defn}

We recall a basic definition.
\begin{defn} Let $Gr(n,V)$ be the Grassmannian of $n$-planes in a vector space $V$. We write $\gamma^n$ to denote the \textbf{tautological bundle of $n$-planes}. That is, each fiber of $\gamma^n$ over a point $S \in Gr(n,V)$ is precisely the subspace $S \subset V$. We will also write $(\gamma^n)^*$ to denote the dual vector bundle, where each fiber over $S$ is the dual space of $S$.
\end{defn}

\begin{lemma}\label{zerolocusbundle} The expected dimension of $I_z$ is $\frac{n^2}{2}(75n + 13)$.
\end{lemma}
\begin{proof} Consider the holomorphic vector bundle $\wedge^2 (\gamma^n)^* \otimes \C^n$ over $Gr(n,D_z)$. Define an algebraic section $\sigma_z$ where $\sigma_z(S) = \theta|_{S\times S}$. This is well defined since $\theta$ is a skew-symmetric tensor and $D$ is corank $n$. Then $I_z$ is precisely the zero locus of the algebraic section $\sigma_z$. The expected dimension of $I_z$ is then $\frac{n^2}{2}(75n + 13)$, since $I_z$ is cut out by ${n}\choose{2}$$n$ equations.
\end{proof}

\begin{defn} Let $M$ and $N$ be two almost complex manifolds with complex structures $J_M$ and $J_N$ respectively. Let $f: M \rightarrow N$ be a smooth mapping. The \textbf{antiholomorphic differential} is defined as
\begin{center}
$\Delbar f = \frac{1}{2}(df + J_N\circ df\circ J_M)$.
\end{center}
\end{defn}
Notice that $\Delbar f$ vanishing identically is equivalent to $df$ commuting with the complex structures; it is precisely the obstruction to $f$ being a holomorphic mapping. In our specific situation, we have
\begin{center}
$\Delbar f = \overline{\partial}_{J_f}f = \frac{1}{2}(df + J_Z\circ df\circ J_f)$
\end{center}
as our antiholomorphic differential, with $J_Z$ the complex structure on the universal space $Z$ and $J_f = J$ the complex structure obtained by pullback via the isomorphism $df$ with $TX$ and $TZ/D$. We are interested in the situation where $\Delbar f$ is injective. A crucial observation thanks to A. Viktorova is that injectivity of $\Delbar f$ is independent of $J_f$, as implied by the following. 

\begin{lemma}
$\Delbar f$ is injective if and only if the immersion $f$ is totally real, i.e., the differential satisfies $df_x(T_xX) \cap J_Z(df_x(T_xX)) = 0$ for all $x \in X$.
\end{lemma}
\begin{proof} Suppose $f$ is a totally real immersion. If $\Delbar f(v) = 0$, by definition $df(v) = -J_Z\circ df\circ J_f(v)$. Therefore $df(v)$ is an element of $df(TX)$ and $J_Z\circ df(TX)$, so by totally realness, $df(v) = 0$. Since $f$ is also an immersion, $v = 0$. Suppose conversely that $\Delbar f$ is injective. Let $V = df(TX) \cap J_Z(df(TX))$. Notice that $V$ is $J_Z$-invariant. Complexifying $TX$, $TZ$, and $df$, we have by $J_Z$-invariance (i.e., $V$ being a complex vector space) that there exists a $+i$-eigenvector $v \in V$. Moreover, $v = df(u)$ for some $u \in TX\otimes\C$. Since $J_Z v = iv$, we have that $J_{TZ/D}\lbrack v\rbrack$ mod $D$ = $i\lbrack v \rbrack$ mod $D$, and therefore, by definition of the pullback complex structure, that $J_f u = iu$ (since $df$ is complex linear with respect to the complexification). Evaluating the complexified $\Delbar f$ on $u$,we have
\begin{equation*}
\begin{split}
\Delbar f(u) &= df(u) + J_Z\circ df \circ J_f(u) \\
&= v + iJ_Z\circ df(u) \\
&= v + iJ_Z(v) \\ 
&= v - v = 0 
\end{split}
\end{equation*}
By injectivity of $\Delbar f$, we have that $u = 0$, so $v = 0$.
\end{proof}

We are now able to fully express the differential relations we are interested in studying.

\begin{example} Let $X$ be a smooth manifold and $Z$ be the universal space. Let $\mathcal{R}_{ac} \subset J^1(X,Z)$ be defined as the intersection $\mathcal{R}_{imm} \cap \mathcal{R}_{real} \cap \mathcal{R}_{D}$. This is the \textbf{almost complex} relation, where elements are tuples $(x,y,L)$ such that $L$ is injective, totally real, and transverse to the universal distribution $D$.
\end{example}

\begin{example} Let $I \subset Gr(n,D)$ be the isotropy locus of the Grassmannian bundle over the universal space $Z$. Let $X$ be an almost complex manifold. Let $\mathcal{R}_{I} \subset \mathcal{R}_{ac}$ denote the $I$-directed mappings where $\Delbar f(T_xX) \in I$ for all $x \in X$.  This is equivalently the set of tuples $(x,y,L)$ where $L$ is injective, transverse to $D$, totally real, and where $\Delbar L =  \frac{1}{2}(L + J_Z\circ L\circ J_L)$ maps $T_xX$ into $I$ for all $x$. Here, $J_L$ denotes the complex structure pulled back by $L$ from $TZ/D$. This is the \textbf{complex} relation.
\end{example}

Observe that we can make sense of $\Delbar L$ since all we need for a pullback complex structure is a bundle isomorphism. By definition, the pullback complex structure $J_f$ satisfies $df\circ J_f$ mod $D = J_Z\circ df$ mod $D$, so by precomposing $\Delbar f$ with $J_f$, we have that $\Delbar f$ maps automatically \textit{into} $D$. Moreover it sends tangent planes to $J_Z$-invariant planes since $\Delbar f$ anticommutes with the complex structures. The directed immersion relation for $\Delbar f$ mapping tangent planes into the isotropy locus $I$ is then well posed.
\\ \\
The reason of interest is due to the following theorem of Demailly and Gaussier \cite{demailly}.
\begin{thm}(Demailly, Gaussier)
For every compact $n$-dimensional \textnormal{integrable} complex manifold $(X,J)$, there exists an embedding $X \hookrightarrow Z^\R$ transverse to $D$, contained in the real part of $Z$, such that
\end{thm}
\begin{enumerate}[(i)]
\item \textit{$J = J_f$ and $\overline{\partial}_{J}f$ is injective};
\item Im$(\overline{\partial}_Jf)$ \textit{is contained in the isotropic locus $I$ of the torsion operator $\theta$ of $D$, the intrinsically defined algebraic locus in the Grassmannian bundle $Gr(n,D) \rightarrow Z$ of complex $n$-dimensional subspaces in $D$ consisting of those subspaces $S$ such that $\theta|_{S\times S} = 0.$}
\end{enumerate}


\textbf{Remark}. The relation $\mathcal{R}_{imm}$ for immersions, $\mathcal{R}_{D}$ for maps transverse to the distribution, and $\mathcal{R}_{real}$ for totally real maps are all \textit{open} relations. This gives the relation $\mathcal{R}_{ac}$ the structure of a smooth subfibration of $J^1(X,Z)$. In contrast, the relation $\mathcal{R}_{I}$ is \textit{closed} and potentially singular.




\section{Formal integrability}\label{formal}
We are now led to studying these differential relations more closely as we are interested in finding mappings into the isotropy locus. Recall from section \ref{prelim} that there is a tower of affine bundles $... \rightarrow J^r(M,N) \rightarrow J^{r-1}(M,N) \rightarrow ... \rightarrow J^0(M,N) = M \times N \rightarrow M$. 

\begin{defn} Let $\sigma$ be a section of the trivial fibration $M \times N \rightarrow M$. Then there is an induced section $j^r\sigma$ of the fibration $J^r(M,N) \rightarrow M$ by taking all derivatives of $\sigma$ up to order $r$. This is the \textbf{r-th jet} of $\sigma$. A section $\eta$ of $J^r(M,N) \rightarrow M$ such that $\eta = j^r\sigma$ for some section $\sigma$ of $M \times N \rightarrow M$ is said to be \textbf{holonomic}.
\end{defn}
\begin{example} We are only interested in $J^1(M,N)$. A holonomic section in coordinates is then of the form $j^1\sigma(x) = (x, \sigma(x), d\sigma_x)$ for $\sigma$ a section of $M \times N \rightarrow M$.
\end{example}

\begin{defn} Let $\mathcal{R} \subset J^r(M,N)$ be a differential relation. A section $\eta$ of $J^r(M,N) \rightarrow M$ whose image lies in $\mathcal{R}$ is said to be a \textbf{formal solution} of $\mathcal{R}$. If $\eta$ is holonomic with image in $\mathcal{R}$, we say it is a \textbf{genuine solution} of $\mathcal{R}.$
\end{defn} 

Every genuine solution is a formal solution. The goal then is to first see if there are any purely homotopical obstructions to having a formal solution, before looking for a genuine solution. For our particular case, we can analyze the Grassmannian bundle directly; i.e., given a map into the universal space $Z$ and a tangential lift into the Grassmannian bundle, we can try and homotope the tangential lift to a map into the subspace $I$.
\\ \\
\textbf{Remark}. If a smooth manifold $X$ admits a genuine solution to $\mathcal{R}_{imm} \cap \mathcal{R}_{D} \subset J^1(X,Z)$ for $Z$ the universal space and $D$ the universal distribution, then $X$ will have an almost complex structure. The genuine solution will correspond precisely to an immersion transverse to $D$, thereby pulling back the complex structure $J_{TZ/D}$ to an almost complex structure on $X$. Conversely, if we already have an almost complex structure, then by Demailly and Gaussier's theorem, we have a genuine solution to $\mathcal{R}_{imm} \cap \mathcal{R}_{D} \cap \mathcal{R}_{real}$. It's worth noting that any obstructions to having a genuine solution here should correspond to known homotopical obstructions for admitting an almost complex structure, but we have not investigated this point. 
\\ \\
To homotope a given tangential lift through a fiberwise homotopy, we must study the relative homotopy groups $\pi_*(Gr(n,N-n),I_z)$ where $I_z$ denotes the fiber over $z$ of the subspace $I$, and $Gr(n,N-n)$ is the Grassmannian of $n$-planes in $(N-n)$-space. For dimension reasons, we only need to consider whether there are obstructions up to $* = 2n$, as the obstructions lie in $H^*(X,\pi_{*-1}(Gr(n,N-n),I_z))$ by classical obstruction theory. However, as $I$ may have possibly singular fibers, the obstruction cocycles may not be a priori well defined. In private communication, Demailly informed the author that $I$ is indeed not smooth. It is then not obvious that the ``bundles'' of interest are amenable to the methods of obstruction theory, as they may not be bundles (or even Serre fibrations).
\\ \\
However, it is here where work of Clemente \cite{clemente} comes into play.
\begin{thm} (Clemente)
There exists a Zariski open subbundle $Gr^o(n,D)$ of the Grassmannian bundle $Gr(n,D)$ over the universal embedding space $Z$ of Demailly and Gaussier such that the statement of Theorem \ref{DemaillyTheorem} holds after replacing $Gr(n,D)$ with $Gr^o(n,D)$. Likewise, the isotropy locus $I$ can be replaced with the bundle $I^o$, whose fibers are the intersection of fibers of $I$ and $Gr^o(n,D)$. The bundle $Gr^o(n,D)$ consists of $n$-planes $S$ in $D$ such that $d\pi|_S$ is injective, where $\pi$ is the projection map $Z \rightarrow \C^{8n}$ from the universal Demailly-Gaussier space to its first euclidean component.  
\end{thm}
Indeed, Clemente shows that the above subspaces are actual bundles. Moreover, she shows that the natural universal map for any almost complex manifold discussed in section \ref{prelim} actually has its antiholomorphic differential map into $Gr^o(n,D)$, and that integrability can be expressed as mapping into $I^o$. She also shows the following:

\begin{thm} (Clemente) The bundles $Gr^o(n,D)$ and $I^o$ have the same fiberwise homotopy type.
\end{thm}

Clemente in fact shows both bundles admit the structure of holomorphic affine bundles over a separate Grassmannian bundle. This result immediately shows that classical obstruction theory can be set up, and that the obstructions vanish in all dimensions. The previous version of this paper did not show that the obstruction theory can be set up properly; moreover, assuming that the obstruction theory could be set up in the original context, the previous paper only showed the vanishing of obstructions up to complex dimension 77. This vanishing followed from algebro-geometric considerations using methods of Sommese; for full disclosure's sake, we leave below the arguments of the old paper as they still illustrate some of the algebro-geometric structure of the original spaces $Gr(n,D)$ and $I$.

\begin{defn} Let $E$ be a holomorphic vector bundle over a complex manifold. We say $E$ is \textbf{globally generated} if there exist global holomorphic sections that span the fiber $E_x$ at every point $x$. 
\end{defn} 

\begin{defn} Let $E$ be a holomorphic vector bundle with fiber $E_x$ over $x$. We write $\mathbb{P}E$ to denote the associated \textbf{projective bundle} where each fiber over $x$ is the projective space $\mathbb{P}(E_x)$ of lines in $E_x$. We write $\gamma_E = \mathcal{O}_{\mathbb{P}E}(-1)$ for the tautological bundle of lines on $\mathbb{P}E$. Similarly, we write its dual as $\gamma_E^* = \mathcal{O}_{\mathbb{P}E}(1)$. We will use superscript $s$ to denote tensor product of the bundle with itself $s$ times; e.g., $(\gamma_E^*)^{\otimes s} = \mathcal{O}_{\mathbb{P}E}(s)$ is the tensor product of $s$ copies of $\gamma_E^*$. 
\end{defn}

We also require the following definition as found in \cite{sommese1}.

\begin{defn} Let $E$ be a holomorphic vector bundle over a connected projective manifold $X$. We say $E$ is \textbf{$k$-ample} if $(\gamma_E^{*})^{\otimes s}$ is globally generated for some $s$ and if the map associated to $H^0(\mathbb{P}E, (\gamma_E^{*})^s)$ has fibers of at most dimension $k$.
\end{defn}

The previous version of the paper used the following theorem of Sommese \cite{sommese3}.
\begin{thm} (Sommese)
Let $E$ be a $k$-ample vector bundle on a compact complex manifold $W$. Assume that $E$ is globally generated by sections and that $B$ is the zero set of a holomorphic section of $E$. Then we have $\pi_j(W,B) = 0$ for $j \leq \textnormal{dim} W - \textnormal{rank} E - k$.
\end{thm}
The theorem above is implied by the main results in \cite{sommese3}; a direct proof can also be found in Sommese and Van de Ven \cite{sommese2} (Remark 3.2.1). The following sequence of lemmas and overall argument is due to incredibly helpful conversations with J. Sheridan and J. Starr.
\\ \\
We will require the following theorem as stated in Ch. III, Theorem 12.11 in Hartshorne \cite{hartshorne} on p 290. We will not define the terms in the statement of the theorem, as that would take us too far afield. We only require it for the subsequent corollary \ref{basechangecor}.

\begin{prop}\label{cohomologybasechange} (Cohomology and base change) Let $f: X \rightarrow Y$ be a projective morphism of noetherian schemes and let $\mathcal{F}$ be a coherent sheaf on $X$, flat over $Y$. Let $y$ be a point of $Y$. Then:

\begin{enumerate}[(a)]
\item if the natural map
\begin{equation*}
\phi^i(y): R^if_*(\mathcal{F}) \otimes k(y) \rightarrow H^i(X_y,\mathcal{F}_y)
\end{equation*}
is surjective, then it is an isomorphism, and the same is true for all $y'$ in a suitable neighborhood of $y$.
\item Assume that $\phi^i(y)$ is surjective. Then the following conditions are equivalent:
\begin{enumerate}[(i)]
\item $\phi^{i-1}(y)$ is also surjective;
\item $R^if_*(\mathcal{F})$ is locally free in a neighborhood of $y$.
\end{enumerate}
\end{enumerate}
\end{prop}

We also need the following sequence of lemmas concerning the Grassmannian $Gr(r,s)$.
\\ \\
For the remainder of this section, let $E = (\gamma^r)^*$ denote the dual of the tautological bundle of $r$-planes in $s$-space over $Gr(r,s)$ unless otherwise stated.

\begin{lemma}\label{alggeo1} The projectivization $\mathbb{P}E$ is the point-plane incidence variety in $Gr(r,s) \times \mathbb{P}(\C^s)$. Then there is a natural projection map $\mathbb{P}E \xrightarrow{p} \mathbb{P}(\C^s)$. The fiber of $p$ over a line $\lambda \in \mathbb{P}(\C^s)$ is then $Gr(r-1,\C^s/\lambda)$.
\end{lemma}
\begin{proof} The first statement is immediate from the definition of $\mathbb{P}E$, as the fiber over the $r$-plane $S$ is exactly the projectivization of $S$.
\end{proof} 

\begin{lemma}\label{alggeo2} There is a bundle isomorphism $p^*\mathcal{O}_{\mathbb{P}(\C^s)}(1) \simeq \mathcal{O}_{\mathbb{P}E}(1)$.
\end{lemma} 
\begin{proof} This again follows from the definition of $p$, since over a point $(\Lambda,\lambda) \in \mathbb{P}E$ where $\Lambda$ is an $r$-plane and $\lambda$ is a line in $\Lambda$, each bundle has as fiber all linear functionals defined on $\C^s$ restricted to $\lambda \subset \Lambda$.
\end{proof} 

We can now state the corollary of proposition \ref{cohomologybasechange}.

\begin{cor}\label{basechangecor} If the groups $H^{q}(Gr(r-1,\C^s/\lambda), \mathcal{O}_{\mathbb{P}E}(1)|_{Gr(r-1,\C^s/\lambda)})$ vanish, then so do the groups $H^p(\mathbb{P}(\C^s), R^{q}p_{*}\mathcal{O}_{\mathbb{P}E}(1))$.
\end{cor}
\begin{proof} This follows from part (a) of \ref{cohomologybasechange} applied to $X = \mathbb{P}E$ and $Y = \mathbb{P}(\C^s)$, with $f = p$ the projection, as in the statement of the theorem. The fiber $X_y$ is then $Gr(r-1,\C^s/\lambda)$ as in lemma \ref{alggeo1}. The sheaf $\mathcal{F}$ is the sheaf of sections of the bundle $\mathcal{O}_{\mathbb{P}E}(1)$.
\end{proof}

\begin{lemma} The bundle $\mathcal{O}_{\mathbb{P}E}(1)|_{Gr(r-1,\C^s/\lambda)}$ over $Gr(r-1, \C^s/\lambda)$ is isomorphic to the trivial line bundle $\mathcal{O}_{Gr(r-1,\C^s/\lambda)}$.
\end{lemma}
\begin{proof}  We have the bundle isomorphism $\mathcal{O}_{\mathbb{P}E}(1) \simeq  p^*\mathcal{O}_{\mathbb{P}(\C^s)}(1)$ from lemma \ref{alggeo2}. We are restricting the bundle to $Gr(r-1,\C^s/\lambda)$, which is exactly the fiber of $p$. That is, it is mapped to a point. The pullback of a bundle over a point is trivial.
\end{proof}

We now need to recall some basic definitions and one famous theorem from algebraic geometry.

\begin{defn} Let $L$ be a holomorphic line bundle over a projective variety $X$. We say $L$ is \textbf{basepoint-free} if the intersection of zero sets of all global sections of $L$ is empty. We say $L$ is \textbf{ample} if there is some integer $r$ such that $L^{\otimes r}$ is basepoint-free and the associated morphism $f: X \rightarrow \mathbb{P}H^0(X, L^{\otimes r})$ is a closed immersion.
\end{defn}

\begin{defn} The \textbf{canonical bundle} $\mathcal{K}$ on a smooth variety $X$ is defined as the top exterior power of the cotangent bundle $T^*X$. The \textbf{anticanonical bundle} $-\mathcal{K}$ is the inverse bundle of $\mathcal{K}$.
\end{defn} 

The following theorem is stated from Ch. III, Remark 7.15 in Hartshorne \cite{hartshorne} on p. 248.

\begin{thm}\label{kodairavanishing} (Kodaira vanishing) Let $X$ be a smooth projective variety of dimension $n$ and $\mathcal{K}_X$ be its canonical bundle. If $L$ is an ample line bundle on $X$, then
\begin{enumerate}[(a)]
\item $H^i(X, L\otimes \mathcal{K}_X) = 0$ for $i > 0$
\item $H^i(X, L^{-1}) = 0$ for $i < n$.
\end{enumerate}
\end{thm}

The following fact is well known. We present an argument due to G. Elencwajg found on Math Stack Exchange \cite{georges}.

\begin{lemma} Let $-K_{Gr(r,s)}$ denote the anticanonical bundle of the Grassmannian $Gr(r,s)$. Then $-K_{Gr(r,s)}$ is ample.
\end{lemma} 
\begin{proof}
Let $X = Gr(r,s)$. Let $TX$ denote the tangent bundle and $T^*X$ denote the cotangent bundle. The anticanonical bundle is then $\wedge^{top} TX$. We have that $TX = $ Hom$(\gamma^r, \gamma^{s-r}_Q)$ where $\gamma_Q^{s-r}$ denotes the tautological quotient bundle of rank $s-r$.  Thus we have that $TX = (\gamma^r)^* \otimes \gamma_Q^{s-r}$. This implies that $\wedge^{top} TX = (\wedge^{top} (\gamma^r)^*)^{\otimes s-r} \otimes (\wedge^{top} \gamma_Q^{s-r})^{\otimes r}$. Consider the short exact sequence
\begin{equation*}
0 \rightarrow \gamma^r \rightarrow \mathcal{O}_X \rightarrow \gamma_Q^{s-r} \rightarrow 0.
\end{equation*} 
From this sequence it follows that $\wedge^{top} \gamma_Q^{s-r} = \wedge^{top} (\gamma^r)^*$. We then obtain that $\wedge^{top} TX = (\wedge^{top} (\gamma^r)^*)^{\otimes s}$. Let $P: Gr(r,s) \hookrightarrow \mathbb{P}(\wedge^r \C^s)$ be the Pl\"{u}cker embedding. The pullback $P^*\mathcal{O}_{\mathbb{P}(\wedge^r\C^s)}(1)$ is precisely $\wedge^{top} (\gamma^r)^*$, so we conclude that $\wedge^{top} (\gamma^r)^*$ is ample. It then follows that $TX$ is ample, being the $s$-th tensor power of $\wedge^{top} (\gamma^r)^*$.
\end{proof}

\begin{cor}\label{COHvanishes} The cohomology groups $ H^{q}(Gr(r-1,\C^s/\lambda), \mathcal{O}_{Gr(r-1,\C^s/\lambda)})$ all vanish for $q > 0$.
\end{cor}
\begin{proof} Let $X = Gr(r-1, \C^s/\lambda)$. We have that the anticanonical $-K_X$ is ample. Using Kodaira vanishing with $L = -K_X$, we obtain $H^{q}(X, \mathcal{O}_X) = H^q(X, L\otimes \mathcal{K}_X)= 0$ for all $q > 0$.
\end{proof} 

\begin{cor}\label{highersheavesvanish} The groups $H^p(\mathbb{P}(\C^s), R^{q}p_{*}\mathcal{O}_{\mathbb{P}E}(1))$ all vanish for $q > 0$.
\end{cor} 
\begin{proof} This follows from corollary \ref{COHvanishes} and corollary \ref{basechangecor}.
\end{proof}

We lastly require the following projection formula, again found in p. 124 of Hartshorne \cite{hartshorne} as an exercise 5.1.d in Ch. II. Again we will not define the terms here, for sake of brevity. Our only concern will be its corollary.

\begin{prop}\label{projectionformula} (Projection formula) If $f: (X,\mathcal{O}_X) \rightarrow (Y, \mathcal{O}_Y)$ is a morphism of ringed spaces, if $\mathcal{F}$ is an $\mathcal{O}_X$-module, and if $\mathcal{E}$ is a locally free $\mathcal{O}_Y$-module of finite rank, then there is a natural isomorphism $f_*(\mathcal{F} \otimes_{\mathcal{O}_X} f^*\mathcal{E}) \simeq f_*(\mathcal{F}) \otimes_{\mathcal{O}_Y} \mathcal{E}$.
\end{prop}

\begin{cor}\label{alggeo3}  Let $p_{*}\mathcal{O}_{\mathbb{P}E}(1)$ denote the pushforward sheaf under the projection $\mathbb{P}E \xrightarrow{p} \mathbb{P}(\C^s)$. Then, $p_{*}\mathcal{O}_{\mathbb{P}E}(1) = \mathcal{O}_{\mathbb{P}(\C^s)}(1)$.
\end{cor} 
\begin{proof} We have that $p_{*}\mathcal{O}_{\mathbb{P}E}(1) = p_{*}(\mathcal{O}_{\mathbb{P}E}\otimes p^{*}\mathcal{O}_{\mathbb{P}(\C^s)}(1))$ by lemma \ref{alggeo2}. But this latter sheaf, by the projection formula, is $p_{*}(\mathcal{O}_{\mathbb{P}E})\otimes \mathcal{O}_{\mathbb{P}(\C^s)}(1)$. But this sheaf is isomorphic to $\mathcal{O}_{\mathbb{P}(\C^s)}(1)$, since the fibers of $p$ are connected, which implies that the sheaf for the trivial bundle $\mathcal{O}_{\mathbb{P}E}$ is pushed forward to the sheaf for the trivial bundle $\mathcal{O}_{\mathbb{P}(\C^s)}$ under $p$. 
\end{proof} 

We now need the following corollary regarding the sheaf cohomology of projective space, phrased more generally as Theorem 5.1 in Ch. III of Hartshorne \cite{hartshorne}, on p. 225. We state it as a proposition specific to our case here so as to avoid excessive terminology.

\begin{prop}\label{cohomologyofprojectivespace} $H^p(\mathbb{P}(\C^s),\mathcal{O}_{\mathbb{P}(\C^s)}(1)) = 0$ for $p > 0$ and $H^0(\mathbb{P}(\C^s),\mathcal{O}_{\mathbb{P}(\C^s)}(1)) \simeq (\C^s)^{*}$.
\end{prop}

\begin{lemma}\label{dimensionEqualsS} The dimension of the space of global sections of $p^*\mathcal{O}_{\mathbb{P}(\C^s)}(1)$ over $\mathbb{P}E$, denoted $h^0(\mathbb{P}E, p^*\mathcal{O}_{\mathbb{P}(\C^s)}(1))$,  is precisely equal to $s$.
\end{lemma} 
\begin{proof} To compute the above dimension, we apply the Leray spectral sequence to the point-plane incidence fibration $\mathbb{P}E \xrightarrow{p} \mathbb{P}(\C^s)$ whose fiber over a line $\lambda$ is $Gr(r-1,\C^s/\lambda)$. On the $E_2$ page we have terms $H^p(\mathbb{P}(\C^s), R^{q}p_{*}\mathcal{O}_{\mathbb{P}E}(1))$ which abut to $H^{p+q}(\mathbb{P}E, \mathcal{O}_{\mathbb{P}E}(1))$. However, by corollary \ref{highersheavesvanish} we have that the groups $H^p(\mathbb{P}(\C^s), R^{q}p_{*}\mathcal{O}_{\mathbb{P}E}(1))$ all vanish for $q > 0$.
\\ \\
For $q = 0$, we have the ordinary pushforward sheaf $p_{*}\mathcal{O}_{\mathbb{P}E}(1)$ instead of the higher derived ones. But $H^p(\mathbb{P}(\C^s), p_*\mathcal{O}_{\mathbb{P}E}(1))$ is isomorphic to $H^p(\mathbb{P}(\C^s),\mathcal{O}_{\mathbb{P}(\C^s)}(1))$ by corollary \ref{alggeo3}. These latter groups all vanish for $p > 0$ and is isomorphic to the dual space $(\C^s)^*$ for $p = 0$, by proposition \ref{cohomologyofprojectivespace}. Thus the $E_2$ page of the spectral sequence only has one nonzero term, at $(p,q) = (0,0)$. We conclude that $h^0(\mathbb{P}E, p^*\mathcal{O}_{\mathbb{P}(\C^s)}(1)) = s$.
\end{proof} 

We conclude with one last definition from algebraic geometry.
\begin{defn} A map $f: X \rightarrow \mathbb{P}(\C^s)$ is \textbf{nondegenerate} if it does not map into a hyperplane.
\end{defn} 
It is well known that a map $f$ to projective space $\mathbb{P}(\C^s)$ is given by taking all global sections of $f^*\mathcal{O}_{\mathbb{P}(\C^s)}(1)$ if and only if $f$ is nondegenerate and the dimension of the space of global sections of $f^*\mathcal{O}_{\mathbb{P}(\C^s)}(1)$ is equal to $s$. See Griffiths \& Harris \cite{griffiths}.

\begin{lemma}\label{globalsectionsmap} The natural projection map $\mathbb{P}E \xrightarrow{p} \mathbb{P}(\C^s)$ is the map associated to global sections of $\mathcal{O}_{\mathbb{P}E}(1)$. That is, $p$ is nondegenerate and $h^0(\mathbb{P}E, p^*\mathcal{O}_{\mathbb{P}(\C^s)}(1)) = s$.
\end{lemma}
\begin{proof} The map $p$ is nondegenerate, since it is the identity in the projective space component, being a projection. By lemma \ref{dimensionEqualsS}, we have $h^0(\mathbb{P}E, p^*\mathcal{O}_{\mathbb{P}(\C^s)}(1)) = s$. Thus it is the associated map to global sections of $\mathcal{O}_{\mathbb{P}E}(1)$.
\end{proof}

\begin{prop}
The bundle $(\gamma^r)^*$ over the Grassmannian $Gr(r,s)$ of $r$-planes in $s$-space is $k$-ample for $k = (r-1)(s-r)$.
\end{prop} 
\begin{proof} Let $E = (\gamma^r)^*$. Its projectivization $\mathbb{P}E$ is the point-plane incidence correspondence in $Gr(r,s) \times \mathbb{P}(\C^s)$, since one is projectivizing each $r$-plane tautologically. Thus there is a natural projection $\mathbb{P}E \xrightarrow{p} \mathbb{P}(\C^s)$. The bundle $\gamma_E^* = \mathcal{O}_{\mathbb{P}E}(1)$ is already globally generated, so we are interested in the map associated to $H^0(\mathbb{P}E, \mathcal{O}_{\mathbb{P}E}(1))$.
\\ \\
By lemma \ref{globalsectionsmap}, the projection $\mathbb{P}E \xrightarrow{p} \mathbb{P}(\C^s)$ is exactly the map associated to global sections of $\mathcal{O}_{\mathbb{P}E}(1)$. The fiber of $p$ has dimension $k = (r-1)(s-r)$.
\end{proof}

We now return from algebraic geometry to our particular case with $I_z \subset Gr(n,D_z)$ over a point $z$ in the universal space $Z$. Recall from lemma \ref{zerolocusbundle} that $I_z$ is the zero locus of an algebraic section of a holomorphic bundle over $Gr(n,D_z)$. This is how the previous paper obtained the 77 dimension bound.

\begin{lemma} The relative homotopy groups $\pi_j(Gr(n,N-n),I_z)$ vanish for $j \leq \frac{1}{2}(-n^3 + 77n^2 + 12n)$.
\end{lemma}
\begin{proof} Let $E = \wedge^2 (\gamma^n)^* \otimes \C^n$ over $Gr(n,D_z)$. Since it is a quotient of the trivial bundle, $(\gamma^n)^*$ is globally generated. Being $k$-ample is closed under the operations of direct sum, quotient, and tensor \cite{lazarsfeld}, for globally generated bundles. So is the property of being globally generated. Thus, our bundle $E$ is globally generated and $k$-ample for $k = 38n^3 - 32n^2 - 6n$. The relative homotopy groups then vanish by Sommese's theorem.
\end{proof}

The following construction is due to a helpful discussion with A. Viktorova, who very helpfully explained how to obtain the bundle automorphisms on $D$.

\begin{lemma}\label{bundleautos} Let $Z$ be the universal space with $D$ its universal distribution. Let $X$ be an almost complex manifold with an immersion $f: X \rightarrow Z$ and tangential lift $\Delbar f$ into $Gr^o(n,D)$. Then there is a family of complex bundle automorphisms $\Phi_t: D \rightarrow D$ that cover the Identity on $Z$ such that $\Phi_0 = Id$ and $\Phi_1\circ\Delbar f$ maps into the isotropy locus $I^o$.
\end{lemma}
\begin{proof} Equip $Z$ with a hermitian metric and consider the fibration $SU(D)$ over $Z$, where a fiber over $z$ is the space $SU(D_z)$. This is a fiber subbundle of Hom$(D,D)$, not to be confused with a bundle of frames; indeed, $SU(D)$ admits sections. Then we can define a projection map $SU(D) \rightarrow Gr^o(n,D)|_{f(X)}$ by the fiberwise quotient map $SU(N-n) \rightarrow SU(N-n) / (SU(n) \times SU(N-2n))$. Note this quotient map is only well defined as a fibration over $Gr(n,D)|_{f(X)}$ since otherwise there are no canonical choices of $n$-planes to take quotients by. With the immersion, we can choose our $n$-planes to be $\Delbar f(T_xX)$.
By Clemente's theorem, we have a fiberwise deformation retract $H: f(X) \times [0,1] \rightarrow Gr^o(n,D)|_{f(X)}$ that deforms the tangential map $\Delbar f$ into the isotropy locus $I^o$.  We then apply homotopy lifting to lift the tangential homotopy $f(X) \times [0,1] \rightarrow Gr(n,D)^o|_{f(X)}$ to a homotopy $f(X) \times [0,1] \rightarrow SU(D)$, by lifting the initial time $H(\cdot, 0) = \Delbar f$ to the Identity map in $SU(D)$. We now have a homotopy $\widetilde{H}: f(X) \times [0,1] \rightarrow SU(D)$ with $\widetilde{H}(\cdot, 0) =$ Id$_D$.
\\ \\
We have that $SU(D) \rightarrow Z$ has a canonical section over $Z$, given by choosing the identity Id$_D$ at each point. By homotopy extension with respect to $(Z, f(X))$, this extends $\widetilde{H}$ to a homotopy $Z \times [0,1] \rightarrow SU(D)$. For each time $t$, we then obtain a section $\widetilde{H}(\cdot, t) = \Phi_t$ of $SU(D)$. This gives us a family of complex bundle automorphisms $\Phi_t: D \rightarrow D$ covering the Identity on $Z$, such that $\Phi_0 = Id$ and $\Phi_1\circ\Delbar f$ maps the tangent planes of $X$ into the isotropy locus.
\end{proof}

\begin{thm}\label{formalsolutionsexist} There are always formal solutions to the differential relation $\mathcal{R}_I$.
\end{thm}
\textit{Proof.} By Clemente's theorem there is a fiberwise deformation retract from $Gr^o(n,D_z)$ to $I^o_z$. Equip $Z$ with a hermitian metric and take the family of bundle automorphisms $\Phi_t$ as constructed in lemma \ref{bundleautos}. Choose a $J_Z$-invariant complementary subbundle of $D$. Extend $\Phi_t$ to be the identity on the complementary subspace and consider $\Phi_t\circ df$. Since we extended as the Identity on $D^\perp$, we have $\Phi_t\circ df$ mod $D = df$ mod $D$, so we also have for the pullback complex structures that $J_{\Phi_t\circ df} = J_f$. Since $\Phi_t$ is $J_Z$-linear, we conclude $\Delbar_{J_{\Phi_t\circ df}}(\Phi_t\circ df) = \Phi_t\circ\Delbar_{J_f} f$. Moreover, $\Phi_t\circ df$ is still an injective, totally real linear map transverse to $D$ for all $t$. Finally, $\Phi_t\circ\Delbar f(TX)$ is still  contained in  $Gr^o(n,D)$ since the family $\Phi_t$ in $SU(D)$ projects exactly to the homotopy in $Gr^o(n,D)$ by the images of $\Phi_t$ on $\Delbar f(TX)$. $\qed$

\begin{defn} Let $\mathcal{R}$ be a relation in $J^1(X,Z)$. Let $\Gamma(X,\mathcal{R})$ denote the space of sections of $J^1(X,Z) \rightarrow X$ whose image is in $\mathcal{R}$.
\end{defn}

\begin{cor} There is a deformation retract from $\Gamma(X, \mathcal{R}_{ac})$ to $\Gamma(X, \mathcal{R}_I)$. 
\end{cor}
\begin{proof} In the argument of lemma \ref{bundleautos}, the constructed family $\Phi_t = \Phi^f_t$ depends continuously on the $0$-jet component $f$ and its tangential lift $\Delbar f$. More generally, it depends continuously on the whole $1$-jet section. Define a homotopy $F: \Gamma(X, \mathcal{R}_{ac}) \times [0,1] \rightarrow \Gamma(X, \mathcal{R}_{ac})$ where $F(\sigma, t) = \Phi^\sigma_t \circ\sigma$, where if $\sigma = (x,z,L)$ then $\Phi^\sigma_t\circ\sigma(x) = (x,z,\Phi^\sigma_t\circ L)$. The claim is that this is a deformation retraction.
\\ \\
The space $\Gamma(X,\mathcal{R}_I)$ clearly includes into $\Gamma(X,\mathcal{R}_{ac})$. For sections $\sigma$ that already map to $\mathcal{R}_I$, the constructed family $\Phi^\sigma_t$ is the Identity for all time $t$. We have that $F(\sigma,0) = \sigma$ since $\Phi_0 = Id$ over all $\sigma$. Moreover, $F(\sigma, 1) = \Phi^\sigma_1\circ\sigma$ which maps into $\Gamma(X,\mathcal{R}_I)$ by construction. Thus $F$ is a homotopy between the retract $F(\cdot, 1)$ and the Identity on $\Gamma(X,\mathcal{R}_{ac})$.
\end{proof} 

We now recall there is a natural action of the diffeomorphism group Diff $X$ on $J^1(X,Z)$. For a $1$-jet $(x,z,L)$, and a diffeomorphism $g$ in Diff $X$, the action is defined as $g\cdot(x,z,L) = (g(x),z, L\circ dg^{-1})$. This action then naturally extends to sections of $J^1(X,Z)$.

\begin{defn} We say a relation $\mathcal{R} \subset J^1(X,Z)$ is \textbf{\textnormal{Diff} $X$-invariant} if $g\cdot \mathcal{R} \subset \mathcal{R}$ for all $g$ in \textnormal{Diff} $X$.
\end{defn}
It is immediate that $\mathcal{R}_{ac}$ and $\mathcal{R}_c$ are Diff $X$-invariant relations.

\begin{lemma} The deformation retract from $\Gamma(X,\mathcal{R}_{ac})$ to $\Gamma(X,\mathcal{R}_I)$ is \textnormal{Diff} $X$-invariant. That is, if $g$ is an element of \textnormal{Diff} $X$, then $g\cdot F(\sigma, t) = F(g\cdot \sigma, t)$ for all time $t$.
\end{lemma}
\begin{proof} Recall that the deformation retract is defined as $F(\sigma, t) = \Phi^\sigma_t\circ\sigma$, where if $\sigma(x) = (x,z,L)$ then $\Phi^\sigma_t\circ\sigma(x) = (x,z,\Phi^\sigma_t \circ L)$. Fix a time $t$. Then $(g\cdot F(\sigma, t))(x) = (g(x), z, \Phi^\sigma_t\circ L \circ dg^{-1})$. On the other hand, we have that $F(g\cdot\sigma, t)(x) = (g(x), z, \Phi^{g\cdot\sigma}_t \circ L \circ dg^{-1})$. The only thing to check is that $\Phi^\sigma_t = \Phi^{g\cdot\sigma}_t$. But this is true since the constructed bundle automorphisms $\Phi^\sigma_t$ only depend on the image of the 0-component of $\sigma$ and the images of the tangent planes under all the $\Delbar L$. But $\sigma$ and $g\cdot\sigma$ have the same 0-component, and the images of the tangent planes under all the $L$ and $L \circ dg^{-1}$ are the same.
\end{proof}

\textbf{Remark.} We conclude that there are always formal solutions to the directed immersion problem. The question remains as to when $\Phi_1\circ df$ is actually a genuine solution, or at least homotopic to a genuine solution. In fact, $\Phi_1\circ df$ as constructed above in theorem \ref{formalsolutionsexist} will \textit{never} be holonomic unless the initial complex structure was integrable to begin with, since the pull back complex structure by the constructed formal solution \textit{is} the initial complex structure. There may be no genuine solutions at all, as mentioned in the introduction. Nonetheless, we can try and \textit{approximate} the formal solution by a holonomic section, and see how close we can get to the isotropy locus.

\section{Holonomic approximation of a complex structure}\label{ample}
In this section, we prove that given a formal solution, we can always find a holonomic section that is $\epsilon$-close, albeit after perturbing our manifold within a thickened neighborhood of itself. If $A$ is a subset, we denote by $Op\ A$ an open neighborhood of $A$, following Gromov's notation \cite{gromov}. We now introduce two technical definitions, which can be found in Eliashberg and Mishachev \cite{eliashberg}.
\begin{defn}
Let $X \rightarrow V$ be a fibration and let $I^k$ denote the $k$-dimensional unit cube. A differential relation $\mathcal{R} \subset X^{(r)}$ is \textbf{locally integrable} if any formal solution over a point can be locally extended to a genuine solution in a parametric and relative sense. That is, given a map $h: I^k \rightarrow V$, a family of sections parameterized by $p$ in $I^k$:
\begin{equation*}
F_p: h(p) \rightarrow \mathcal{R}, \ p \in I^k
\end{equation*}
and a family of local holonomic extensions near $\partial I^k$:
\begin{equation*}
\tilde{F}_p: Op\ h(p) \rightarrow \mathcal{R},
\end{equation*}
\begin{equation*}
\tilde{F}_p(h(p)) = F_p(h(p)), \ p \in Op(\partial I^k),
\end{equation*}
there exists a family of local holonomic extensions
\begin{equation*}
\tilde{F}_p: Op\ h(p) \rightarrow \mathcal{R},
\end{equation*}  
\begin{equation*}
\tilde{F}_p(h(p)) = F_p(h(p)), \ p \in I^k
\end{equation*}
such that for $p \in Op(\partial I^k)$, these new extensions agree with the original extensions over $Op\ h(p)$. 
\end{defn}
\begin{defn} Fix $n =$ dim $V$. Let a $\theta_k$-\textbf{pair} be any pair $(A,B)$ diffeomorphic to $(\lbrack -1, 1\rbrack^n, \lbrack -1, 1\rbrack^k \cup \partial(\lbrack -1, 1 \rbrack^n))$.
\end{defn}
\begin{defn}
A differential relation $\mathcal{R}$ is \textbf{microflexible} if local deformations of genuine solutions can be extended to global deformations of genuine solutions for small times. That is, for any sufficiently small open ball $U \subset V$ and any families smoothly parameterized by $p \in I^m$ of
\begin{itemize}
\item $\theta_k$-pairs $(A_p, B_p) \subset U$,
\item holonomic sections $F^0_p: Op \ A_p \rightarrow \mathcal{R}$, and
\item holonomic homotopies $F^\tau_p: Op \ B_p \rightarrow \mathcal{R}, \tau \in \lbrack 0, 1\rbrack$, of the sections $F^0_p$ over $Op \ B_p$ which are constant over $Op(\partial B_p)$ for all $p \in I^m$ and constant over $Op \ B$ for $p \in Op(\partial I^m)$,
\end{itemize}
 there exists a number $\sigma > 0$ and a family of holonomic homotopies
\begin{equation*}
F^\tau_p: Op \ A_p \rightarrow \mathcal{R}, \ \tau \in \lbrack 0, \sigma\rbrack,
\end{equation*}
which extend the family of homotopies
\begin{equation*}
F^\tau_p: Op \ B_p \rightarrow \mathcal{R}, \ \tau \in \lbrack 0, \sigma\rbrack,
\end{equation*}
and are constant over $Op(\partial A_p)$ for all $p \in I^m$ and constant over $Op \ A$ for $p \in Op(\partial I^m)$.
\end{defn}
Microflexibility will be used to glue/interpolate from one holonomic section to another when both are defined on the same open neighborhood.
\\ \\
We can now state the following theorem also in Eliashberg and Mishashev \cite{eliashberg}. 
\begin{thm} (Holonomic $\mathcal{R}$-approximation theorem) Let $\mathcal{R} \subset X^{(r)}$ be a locally integrable microflexible differential relation. Let $A \subset V$ be a polyhedron of positive codimension and suppose there is a section $F: Op \ A \rightarrow \mathcal{R}$. Assume $V$ is equipped with a metric and the bundle $X^{(r)}$ admits a Euclidean structure in a neighborhood of the section $F(V) \subset X^{(r)}$. Then for arbitrarily small $\delta, \epsilon > 0$, there exists a $\delta$-small (in the $C^0$-sense) diffeotopy $h^\tau: V \rightarrow V, \tau \in \lbrack 0, 1\rbrack$, and a holonomic section $\tilde{F}: Op \ h^1(A) \rightarrow \mathcal{R}$ such that
\begin{equation*}
||\tilde{F}(v) - F|_{Op\ h^1 (A)}(v)||_{C^0} < \epsilon
\end{equation*}
for all $v \in Op \ h^1(A)$.
\end{thm} 
Given a formal solution to a locally integrable microflexible differential relation defined on some positive codimension polyhedron, we can find an $C^0$-approximating holonomic solution \textit{in the $r$-th jet space $X^{(r)}$} defined on a \textit{perturbation} of the polyhedron. In fact, in the case of $V = A \times \R$, we can choose our diffeotopy to be a vertical perturbation (i.e., of the form $(x,t) \mapsto (x, h(x,t))$. The diffeotopy itself will only be $C^0$-small. In our situation, we will have $V = X \times \R$ for $X$ almost complex equipped with metric $g$, and $V$ equipped with the metric $g_{ij} \oplus dt^2$. Equipping the universal Demailly-Gaussier space $Z$ with a metric, we obtain a canonical metric on all jet spaces $J^r(X \times \R, Z)$.  We will need the following exercise of Gromov, found in p. 84 of \cite{gromov}.
\begin{exercise} (Gromov) Let $\mathcal{R}_{tang}$ denote the differential relation of immersions $\ \R \rightarrow Z$ that are \textnormal{tangent} to the distribution $D$. If $D$ is bracket-generating, then $\mathcal{R}_{tang}$ is microflexible. 
\end{exercise}
The exercise above is actually not stated correctly; Bryant and Hsu show in \cite{bryant} that there exist examples of bracket-generating distributions (e.g., Engel structures) which possess \textit{rigid integral curves}, i.e., curves tangent to $D$ that cannot be deformed relative to their ends. We thus need the following preliminary definitions.

\begin{defn}
Let $M$ be a smooth manifold equipped with a distribution $D$. Fix a point $p \in M$. We write $C^\infty([0,1],M,D,p)$ for the space of integral maps $\gamma: [0,1] \to (M,D)$ with given initial point $\gamma(0) = p$, endowed with the $C^\infty$-topology.

We define the \textbf{endpoint map}
\begin{align*}
\Endpoint_p: C^\infty([0,1],M,D,p) \quad& \longrightarrow\quad (M,\xi) \\
\Endpoint_p(\gamma) \quad& := \quad \gamma(1).
\end{align*}
\end{defn}
In fact, the space $C^\infty([0,1],M,D,p)$ is a Fr\'echet manifold locally modelled on $C^\infty(\gamma,\gamma^*D)$, where $\gamma$ is any integral curve, and the endpoint map is smooth. This motivates us to look at the following definition:
\begin{defn}
A curve $\gamma \in C^\infty([0,1],M,D,p)$ is \textbf{regular} if the endpoint map $\Endpoint_p$ is a submersion at $\gamma$. Otherwise, $\gamma$ is said to be \textbf{singular}.
\end{defn}
The condition for a curve to be regular cannot be expressed as a relation in a finite jet space. Rather, it is a relation in $J^\infty([0,1],M)$ the space of germs of curves (suitably topologized). A curve being singular is equivalent to its lift to the cotangent bundle being contracted with the canonical symplectic form to zero (see \cite{pino}). This yields a set of \textbf{singularity equations} that the curve must satisfy. What we have are the following:
\begin{defn} Let $\mathcal{R}_{sing-tang}$ in $J^\infty(\R,Z)$ denote the differential relation of immersions $\R \rightarrow Z$ that are \textnormal{tangent} to a distribution $D$ and \textnormal{singular}. Let $\mathcal{R}_{sing-tang,r}$ denote the differential relation in $J^r(\R,Z)$ that are $r$-jets of curves satisfying the singularity equations up to order $r$. Let $\mathcal{R}_{reg-tang}$ and $\mathcal{R}_{reg-tang,r}$ denote the complements respectively. 
\end{defn}

We reformulate Gromov's proposition as follows.
\begin{thm}\label{alvaro} (del Pino, S.)
If $D$ is bracket-generating and real analytic, then for sufficiently large $r$, the relation $\mathcal{R}_{reg-tang,r}$ in $J^r(\R,Z)$ is locally integrable and microflexible.
\end{thm}
We note that one does not need \textit{bracket-generating} as an assumption on $D$ to prove local integrability when considering \textit{all} integral curves; the assumption is used when restricting to \textit{regular} curves. The authors have also only proved the above lemma for the case of analytic $D$; the smooth case remains open. Given the above result, combined with the fact that our distribution $D$ is bracket-generating, we can now prove the main theorem of the paper, following the same argument as in \cite{eliashberg}, \cite{gromov}, \cite{pino}. We prove some preliminary lemmas first.



\begin{defn} Let $X$ be an almost complex manifold. Consider $X \times \R$ and the modified differential relation in $J^\infty(X \times \R, Z)$: we define $\mathcal{R}_{c-tang}$ as immersions $X \times \R \rightarrow Z$ transverse to $D$, but totally real in the horizontal direction, and \textit{tangent} to $D$ and \textit{regular} when restricted to each fiber $\{x\} \times \R$. Likewise, define $\mathcal{R}_{c-tang,r}$ as the modified differential relation in $J^r(X\times\R,Z)$ of immersions $X \times \R \rightarrow Z$ transverse to $D$, totally real in the horizontal direction, and such that the map is a section of $\mathcal{R}_{reg-tang,r}$ when restricted to each fiber $\{x\} \times \R$. Let $\mathcal{R}_{c,r}$ denote the preimage of $\mathcal{R}_c$ under the projection map $J^r(X,Z) \rightarrow J^1(X,Z)$. 
\end{defn}
The main point of the argument is that holonomic approximation applies to $\mathcal{R}_{c-tang}$ after passing to a large enough jet space, since the conditions in the horizontal directions are all open.

In the following definition and lemmas, we will fix $r$ sufficiently large. 


\begin{lemma}\label{singlepointlemmaLocalInt} For $r$ sufficiently large, $\mathcal{R}_{c-tang,r}$ is locally integrable when the parameter is a point. 
\end{lemma}
\begin{proof}
Since $D$ is bracket-generating, it follows from theorem \ref{alvaro} that for $r$ sufficiently large, $\mathcal{R}_{reg-tang,r}$ in $J^r(\R, Z)$ is both locally integrable and microflexible in the full parametric and relative sense. The openness of $\mathcal{R}_{c,r}$ and local integrability of $\mathcal{R}_{reg-tang,r}$ imply the local integrability of ${\mathcal{R}}_{c-tang,r}$: we will prove first the case when the parameter is a single point.
\\ \\
For a formal solution $F = ((x,t),z,L(x)\oplus V(t), \textit{etc})$ of $\mathcal{R}_{c-tang,r}$ over a point $(x,t)$, we can find a small contractible neighborhood around $x$ and a holonomic section $(x, f(x), df_x, \textit{etc})$ of $\mathcal{R}_{c,r}$ agreeing with $(x,z,L(x), \textit{etc horizontal})$ over $x$, since open relations are locally integrable. Then we can use the local integrability of $\mathcal{R}_{reg-tang,r}$ to find a locally defined regular integral curve $\gamma$ that agrees with $V$ and $etc \ vertical$ over $t$, nowhere vanishing on the image of $f$, which is contractible. Our holonomic section agreeing with $F$ over $(x,t)$ is then defined locally as $((x,t), f(x), df_x\oplus\gamma'(t), etc)$.
\end{proof}

\begin{lemma} For $r$ sufficiently large, $\mathcal{R}_{c-tang,r}$ is locally integrable in the full parametric and relative sense.
\end{lemma}
\begin{proof} Take any family of sections $F_p$ of $\mathcal{R}_{c-tang,r}$ parameterized by $h: I^k \rightarrow X \times \R$ for $p$ in $I^k$ with given holonomic extensions $\widetilde{F}_p$ defined on neighborhoods $Op \ h(p)$ for $p$ in $Op(\partial I^k)$. Restrict the $F_p$ to the horizontal direction as in lemma \ref{singlepointlemmaLocalInt}. Since $\mathcal{R}_{c,r}$ is locally integrable in the full parametric and relative sense, we can extend the sections in the horizontal direction such that they agree with the given horizontal extensions for $p$ in $Op(\partial I^k)$. Since $\mathcal{R}_{reg-tang,r}$ is also locally integrable in the full parametric and relative sense, we can then find vertical sections that extend the vertical components of the given family, and agree with the vertical components for $p$ in $Op(\partial I^k)$. The extensions $\widetilde{F}_p$ are then defined as the sum of the horizontal and vertical components, as in lemma \ref{singlepointlemmaLocalInt}. Since the components agree with the horizontal and vertical components of the sections for $p$ in $Op(\partial I^k)$ respectively, the newly defined extensions themselves agree with the given extensions. This proves $\mathcal{R}_{c-tang,r}$ is locally integrable in the full parametric and relative sense.
\end{proof}

We now proceed to prove microflexibility for $\mathcal{R}_{c-tang,r}$. The idea is that one uses microflexibility of $\mathcal{R}_{reg-tang,r}$ to obtain a family of curves over $A$ such that, when glued together, is the desired holonomic homotopy (i.e., we foliate the desired holonomic homotopy with our curves); openness of $\mathcal{R}_{c,r}$ guarantees the sections remain holonomic in the transverse direction. Again we will prove first the case when the parameter is a point, and then generalize to the full parametric and relative case.

\begin{lemma}\label{thetapair}
Suppose we have a $\theta_k$-pair $(A,B) \subset U \subset X \times \R$ where $U$ is a sufficiently small open ball, a holonomic section $F^0$ defined on $Op \ A$, and a holonomic homotopy $F^\tau$ defined on $Op \ B$ and constant on $Op(\partial B)$. Then there is a smooth family $\widetilde{F}^\tau_{x}$ of regular integral homotopies defined on $(\{x\} \times \R) \cap Op \ A$, for $x$ such that $(x,t)$ is in $A$, and which agree with $F^\tau$ when restricted to $(\{x\} \times \R) \cap Op \ B$, up to a uniform small time $\sigma$.
\end{lemma}
\begin{proof}
As the holonomic homotopy $F^\tau$ is defined on an open neighborhood of $B$, we may assume without loss of generality that $B$ meets each fiber $\{x\} \times \R$ transversely (in the argument, we will parameterize a family of curves by $B$; without transversality, we simply parameterize by a face transverse to $B$ within the neighborhood $Op \ B$). By taking a sufficiently fine subdivision into pairs where each face is transverse to the fibers, we may also assume that $B$ intersects each fiber transversely in at most one point.
\\ \\
We restrict $F^\tau$ to $\{x\} \times \R$ for points $(x,t)$ on $A$. This gives us a holonomic homotopy $F|_{\{x\} \times \R}^\tau$ over midpoints in the interior of curves $(Op (\{x\} \times I_x)) \cap B$ in $\mathcal{R}_{reg-tang,r}$, where $I_x$ denotes a small interval over $x$ in $Op \ A$. Treating the curves $I_x$ and their midpoints as theta pairs in each fiber, we have a collection of theta pairs smoothly parameterized by those $x \in X$ such that $\{x\} \times \R$ intersects $B$ nontrivially (so the parameter space is diffeomorphic to $B$), with holonomic homotopies along the midpoints of each theta pair.
\\ \\
By microflexibility of $\mathcal{R}_{reg-tang,r}$ applied to these theta pairs parameterized by $x$, these holonomic homotopies $F|_{\{x\} \times \R}^\tau$ extend to $\widetilde{F}|_{\{x\}\times\R}^\tau$ for a small time $\sigma$ uniform over the parameter $x$, over the curves $(\{x\} \times I_x)$. This extension agrees with the homotopy on the midpoints $(\{x\} \times \R) \cap Op \ B$ and agrees with $F^0$ on $Op(\partial B)$. Moreover, we can extend them to the rest of $(\{x\} \times \R) \cap Op \ A$ by defining them to be $F^0$ restricted to the fibers. This gives a smooth family $\widetilde{F}^\tau_{x}$ of regular integral homotopies $\widetilde{F}^\tau_{x}(t) = \widetilde{F}|_{\{x\}\times\R}^\tau(t)$ parameterized by $x \in B$ and defined on $(\{x\} \times \R) \cap Op \ A$, and which agree with $F^\tau$ when restricted to $(\{x\} \times \R) \cap Op \ B$, up to a uniform small time $\sigma$.
\end{proof}

\begin{lemma}\label{micropoint}
For $r$ sufficiently large, the relation $\mathcal{R}_{c-tang,r}$ is microflexible when the parameter is a point.
\end{lemma}
\begin{proof} Define the desired holonomic homotopy $\widetilde{F}^\tau$ of $\mathcal{R}_{c-tang,r}$ as the smooth family defined in lemma \ref{thetapair}. That is, $\widetilde{F}^\tau(x, t) = \widetilde{F}^\tau_{x}(t)$. The uniform time $\sigma$ is then the one obtained by microflexibility of $\mathcal{R}_{reg-tang,r}$. Since the extended homotopies by the curves are small perturbations interpolating between $F^0$ and $F^\tau$, we have that $\widetilde{F}^\tau$ restricted to the horizontal direction remains in $\mathcal{R}_{c,r}$, since this relation is \textit{open}. Therefore $\mathcal{R}_{c-tang,r}$ is microflexible when the parameter is a point.
\end{proof}

\begin{lemma} For $r$ sufficiently large, the relation $\mathcal{R}_{c-tang,r}$ is microflexible in the full parametric and relative sense.
\end{lemma}
 \begin{proof} Given smooth families of $\theta_k$-pairs $(A_p,B_p)$, holonomic sections $F^0_p: A_p \rightarrow \widetilde{\mathcal{R}}_{c-tang}$, and holonomic homotopies $F^\tau_p$ defined on $Op \ B_p$ for $p$ in $I^m$, we simply reapply the arguments in lemma \ref{thetapair} and lemma \ref{micropoint}. We restrict $F^\tau_p$ to $\{x\} \times \R$ for $(x,t)$ on $A_p$ for all $p$, and again obtain a smooth family of holonomic homotopies over $B_p$ and parameterized smoothly in $p$. Since $\mathcal{R}_{reg-tang,r}$ is microflexible in the full parametric and relative sense, there is an extension for small time which is uniform over all $x$ and $p$. Moreover, these extensions are constant over $Op(\partial A_p)$ for all $p$ in $I^m$ and constant over $Op \ A$ for $p \in Op(\partial I^m)$. The desired extensions are then the smooth families over $x$, for each $p$. Thus $\mathcal{R}_{c-tang,r}$ is microflexible in the full parametric and relative sense. 
 \end{proof}
We can now state the correct theorem proven by the original paper. 
\begin{thm} Let $X$ be any almost complex manifold with metric $g$. For any formal solution $\sigma$ of $\mathcal{R}_I$ and any $\epsilon > 0$, there exists a smooth function $f: X \rightarrow \R$ such that on the graph $\Gamma_f$ of $f$ in $X \times \R$, there exists a holonomic section $\hat{\sigma}$ on $\Gamma_f$ of $\pi^*\mathcal{R}_I$ such that $||\pi^*\sigma - \hat{\sigma}||_{C^0} < \epsilon$ with norm taken in the graph metric with respect to $f$, and $\pi: \Gamma_f \rightarrow X$ the projection from the graph.
\end{thm}
\noindent \textit{Proof.}
Take a formal solution of $\mathcal{R}_{I} \subset \mathcal{R}_{c}$ on $X$. Let $\mathcal{R}_{I,r} \subset \mathcal{R}_{c,r}$ denote the preimage of $\mathcal{R}_I$ for some sufficiently large $r$. We can define a formal solution of $\mathcal{R}_{I,r} \subset \mathcal{R}_{c,r}$ since the fiber of the preimage is contractible (there are no conditions on higher derivatives). We take a triangulation of $X$ and proceed by induction over the dimensions of the simplices.  For the base case of 0-simplices $\Delta^0$, extend the formal solution of $\mathcal{R}_{I,r} \subset \mathcal{R}_{c,r}$ to a formal solution of $\mathcal{R}_{c-tang,r}$ on $\Delta^0 \times \R$, using that $\Delta^0$ is contractible. By local integrability, we can then find a holonomic section of $\mathcal{R}_{c-tang,r}$ over each 0-simplex.
\\ \\
For the inductive step on $k$-simplices $\Delta^k$, assume we have a holonomic section defined on $\partial\Delta^k$, which is a union of $(k-1)$-simplices. By taking a nowhere vanishing vector field defined on $\Delta^k$ tangent to $D$ as above, we extend the formal solution of $\mathcal{R}_{I,r} \subset \mathcal{R}_{c,r}$ to a formal solution of $\mathcal{R}_{c-tang,r}$ on $\Delta^k \times \R$, with the extension agreeing with the inductively defined holonomic extension along the boundaries of the simplices. Then by the above holonomic approximation theorem, we can perturb an open neighborhood of $\Delta^k \times 0$ by a vertical diffeotopy $h$ and obtain a holonomic section defined on an open neighborhood of $h(\Delta^k \times 0)$ that is $\epsilon$-close to our formal solution in $\mathcal{R}_{I,r}$, again agreeing with the inductively defined holonomic section on the boundary. Restrict the section to $h(\Delta^k \times 0)$ for a holonomic section of $\mathcal{R}_{c,r}$ which is $\epsilon$-close to the original formal solution of $\mathcal{R}_{I,r}$. Moreover, the section is well defined globally on the $k$-skeleton of $X$, as the holonomic extensions agree with the holonomic sections defined on the lower dimensional skeleta. Repeat the argument on $\Delta^{k+1}$. Notice that the vertical diffeotopy $h$ is constructed at each step, along each simplex. Finally, projecting to the first jet space, we end with a holonomic section defined on $h(X \times 0)$ that is $\epsilon$-close to the isotropy locus. Since $h$ is a vertical diffeotopy, projecting to $\R$ gives the desired function $f$ for which the statement holds. $\qed$

We reformulate the above in terms of the main theorem \ref{MainTheorem} restated below. 
\begin{cor}
Let $(X, J, g)$ be an almost complex manifold with metric. For any $\epsilon > 0$, there exists a smooth function $f: X \rightarrow \R$ and almost complex structure $J'$ on $X$ such that
\begin{itemize}
\item $||\pi^*J - \pi^*J'||_{C^0} < K\epsilon$ where $\pi^*J$ and $\pi^*J'$ denote the pullbacks of $J$ and $J'$ on $\Gamma_f$ respectively, and where the norm is taken with respect to $\hat{g}_f$ and the constant $K$ only depends on $J$ and $g$
\item the Nijenhuis tensor of $J'$ on $\Gamma_f$ has $C^0$ norm less than $C\epsilon$, again with respect to $\hat{g}_f$, and where $C$ is a constant depending only on $J$ and $g$.
\end{itemize}
\end{cor}
\begin{proof} The constructed holonomic section $\hat{\sigma}$ pulls back the universal complex structure to an almost complex structure on $h(X \times 0)$. Push this almost complex structure forward by $\pi$ to obtain $J'$; the almost complex structure pulled back by $\hat{\sigma}$ is then $\pi^*J'$. The sections $\pi^*\sigma$, $\hat{\sigma}$ of the jet space on $h(X \times 0)$ being $C^0$-close means that the underlying $1$-jets are $C^0$-close. In particular, by taking a trivializing open cover and using compactness to reduce to finite subcovers, we can see in coordinates that the almost complex structures pulled back by the respective sections are $C^0$-close, up to possibly a constant depending only on $J$. This proves the first item. The second item follows from the fact that the section $\hat{\sigma}$ is holonomic and $C^0$ close to the isotropy locus.
\end{proof}
In words, this says: given an almost complex manifold, we can deform the almost complex structure in a particular way so as to make it close to integrable. This is not with respect to the original metric on the manifold, but with respect to the graph metric associated to the deformation. The previous version of this paper's main error conflated these two metrics, which led to the original paper's claim that any almost complex manifold could admit a sequence of almost complex structures whose Nijenhuis norms get arbitrarily close to 0.
\\ \\
In fact, shortly after the previous version of this paper was first announced, Luis Fernandez and Scott Wilson constructed explicit examples of almost complex manifolds that are known to have no integrable structure, with almost complex structures that get arbitrarily close to integrability \cite{fernandez}. We end this paper with a question.
\begin{question} Does every almost complex manifold, with any fixed metric, admit almost complex structures whose Nijenhuis tensors become arbitrarily close to $0$ in the $C^0$ norm?
\end{question} 
One potential strategy to this question may be using convex integration techniques instead of holonomic approximation; an appropriate version of convex integration would deform the map pulling back the universal complex structure in such a way that it may still be $C^1$-close to a formal solution, all pulled back on the original manifold, at the expense of $C^2$. The method of holonomic approximation presented here can only deform the map to be $C^0$-close to the formal solution, on the original manifold, at the unfortunate cost of $C^1$.

\textsc{Stony Brook University, Department of Mathematics}
\\
\textit{E-mail address:} {\fontfamily{cmtt}\selectfont tobias.shin@stonybrook.edu}


\begin{thebibliography}{50}
\bibitem{bryant} Bryant, R. L., and Hsu, L., 1993. \textit{Rigidity of integral curves of rank 2 distributions}, Invent. Math. \textbf{114}, pp. 435-461.
\bibitem{calabi} Calabi, E., and Eckmann, B., 1953. \textit{A class of compact, complex manifolds which are not algebraic}, Ann. Math. \textbf{58}, No. 3, pp. 494-500.
\bibitem{clemente} Clemente, G., 2020. \textit{Geometry of universal embedding spaces for almost complex manifolds}, arXiv preprint, arXiv:1905.06016. 
\bibitem{demailly} Demailly, J.-P., and Gaussier, H., 2017. \textit{Algebraic embeddings of smooth almost complex structures}, J. Eur. Math. Soc. \textbf{19},  pp. 3391-3419. {\url{https://www-fourier.ujf-grenoble.fr/~demailly/documents.html}}[D43] presented at CUNY Graduate Center for Jim Simons's 80th birthday conference. \textit{Accessed:} 2019-04-17.


\bibitem{eliashberg} Eliashberg, Y., and Mishachev, N., 2002. \textit{Introduction to the h-principle}, Grad. Stud. in Math., Vol. 48, American Mathematical Society, Providence, R.I., pp. 44-46, 66-67, 129-133, 136-138.
\bibitem{evans} Evans, J. D., 2012. \textit{The infimum of the Nijenhuis energy}, Math. Res. Lett. \textbf{19}, pp. 383-388.
\bibitem{fernandez} Fernandez, L., Shin, T., Wilson, S., 2021. \textit{Almost complex manifolds with small Nijenhuis tensor}, arXiv preprint, arXiv:2103.06090

\bibitem{georges} Georges Elencwajg. {\url{https://math.stackexchange.com/users/3217/georges-elencwajg}}, Grassmannian a Fano manifold?,  URL (version: 2015-03-24): \url{https://math.stackexchange.com/q/1203337}


\bibitem{glover} Glover, H. H., Homer, W. D., Stong, R. E., 1982. \textit{Splitting the tangent bundle of projective space}, Indiana Univ. Math. J. \textbf{31}, No. 2, pp. 161-166.
\bibitem{goertsches} Goertsches, O., and Konstantis, P., 2017. \textit{Almost complex structures on connected sums of complex projective spaces}, arXiv preprint, arXiv:1710.05316.
\bibitem{griffiths} Griffiths, P., and Harris, J., 1994. \textit{Principles of algebraic geometry}, John Wiley \& Sons, Inc., Hoboken, N.J., pp. 137, 176.
\bibitem{gromov} Gromov, M., 1986. \textit{Partial differential relations}, Ergebnisse der Mathematik und ihrer Grenzgebiete, Vol. 9, Springer-Verlag, Berlin, Heidelberg.
\bibitem{hartshorne} Hartshorne, R., 1977. \textit{Algebraic geometry}, Springer, New York, N.Y., pp. 124, 225, 248, 290-291.
\bibitem{hirsch} Hirsch, M. W., 1976. \textit{Differential topology}, Springer-Verlag, New York, N.Y., pp. 58-66.
\bibitem{hsu} Hsu, L., 1992. \textit{Calculus of variations via the Griffiths formalism}, J. Differential Geometry. \textbf{36}, pp. 551-589.
\bibitem{lazarsfeld} Lazarsfeld, R., 2004. \textit{Positivity in algebraic geometry II. Positivity for vector bundles, and multiplier ideals}, Springer-Verlag, Berlin, Heidelberg, pp. 26-27.
\bibitem{lutian} Lu, P., Tian, and G., 1994. \textit{The complex structure on a connected sum of $S^3 \times S^3$ with trivial canonical bundle}, Math. Ann. \textbf{298}, pp. 761-764.
\bibitem{le} Le, H. V., and Wang, G., 2001. \textit{Anti-complexified Ricci flow on compact symplectic manifolds}, J. Reiene Angew. Math. \textbf{531}, pp. 17-31.
\bibitem{muller} M\"{u}ller, S., and Geiges, H., 2000. \textit{Almost complex structures on 8-manifolds}, Enseign. Math. \textbf{46}, pp. 95-107.
\bibitem{pino}del Pino, \'{A}., and Shin, T., 2020. \textit{Microflexibility and local integrability of horizontal curves}, arXiv preprint, arXiv:2009.14518.
\bibitem{sommese3} Sommese, A. J., 1983. \textit{A convexity theorem}, Proceedings of Symposia in Pure Math. \textbf{40}, Part 2, pp. 497-505.

\bibitem{sommese1} Sommese, A. J., 1978. \textit{Submanifolds of abelian varieties to Rebecca}, Math. Ann. \textbf{233}, pp. 229-256.
\bibitem{sommese2} Sommese, A. J., and Van de Ven, A., 1986. \textit{Homotopy groups of pullbacks of varieties}, Nagoya Math. J. \textbf{102}, pp. 79-90.
\bibitem{van} Van de Ven, A., 1966. \textit{On the Chern numbers of certain complex and almost complex manifolds}, Proc. Natl. Acad. Sci. USA. \textbf{55}, No. 6, pp. 1624-1627.
\bibitem{varolin} Varolin, D., 2001. \textit{The density property for complex manifolds and geometric structures}, J. Geom. Anal. \textbf{11}, No. 1, pp. 135-160.
\bibitem{yang} Yang, H., 2012. \textit{Almost complex structures on (n-1)-connected 2n-manifolds}, Topol. Its Appl. \textbf{159}, pp. 1361-1368.
\end{thebibliography}
\end{document}